\documentclass[thmsa,12pt,sw20lart]{article}%
\usepackage{amsmath}
\usepackage{amsfonts}
\usepackage{amssymb}
\usepackage{graphicx}%
\providecommand{\U}[1]{\protect\rule{.1in}{.1in}}
\newtheorem{theorem}{Theorem}

\newtheorem{corollary}[theorem]{Corollary}

\newtheorem{lemma}[theorem]{Lemma}

\newtheorem{remark}[theorem]{Remark}

\newenvironment{proof}[1][Proof]{\noindent\textbf{#1.} }{\ \rule{0.5em}{0.5em}}
\begin{document}

\author{Dong-Ho Tsai}
\title{Prescribing oscillation behavior of solutions to the heat equation on
$\mathbb{R}^{n}$ via the initial data and its average
integral\thanks{Mathematics\ Subject Classification:\ 35K05, 35K15. }}
\maketitle

\begin{abstract}
Motivated by a classical stabilization result for solution to the Cauchy
problem of the heat equation$\ \partial_{t}u=\bigtriangleup u\ $on
$\mathbb{R}^{n}$, we consider its oscillation behavior with radial initial
data $\varphi\left(  x\right)  =\varphi\left(  \left\vert x\right\vert
\right)  \in C^{0}\left(  \mathbb{R}^{n}\right)  \bigcap L^{\infty}\left(
\mathbb{R}^{n}\right)  .\ $Given four arbitrary finite numbers $r<\alpha
<\beta<s,$ one can construct a radial $\varphi\in C^{0}\left(  \mathbb{R}%
^{n}\right)  \bigcap L^{\infty}\left(  \mathbb{R}^{n}\right)  $ so that
$\varphi\ $together with its corresponding solution$\ u\left(  x,t\right)  $
satisfy the oscillation behavior:
\begin{align*}
\liminf_{\tau\rightarrow\infty}\varphi\left(  \tau\right)   &  =r<\liminf
_{t\rightarrow\infty}u\left(  0,t\right)  =\alpha\\
&  <\limsup_{t\rightarrow\infty}u\left(  0,t\right)  =\beta<\limsup
_{\tau\rightarrow\infty}\varphi\left(  \tau\right)  =s.
\end{align*}
Similarly, given $p<\alpha<\beta<q\ $with $p+q=\alpha+\beta,$ one can find a
radial $\varphi\in C^{0}\left(  \mathbb{R}^{n}\right)  \bigcap L^{\infty
}\left(  \mathbb{R}^{n}\right)  $ so that its average integral $H\left(
\tau\right)  \ $and $u\left(  x,t\right)  $ satisfy%
\begin{align*}
\liminf_{\tau\rightarrow\infty}H\left(  \tau\right)   &  =p<\liminf
_{t\rightarrow\infty}u\left(  0,t\right)  =\alpha\\
&  <\limsup_{t\rightarrow\infty}u\left(  0,t\right)  =\beta<\limsup
_{\tau\rightarrow\infty}H\left(  \tau\right)  =q.
\end{align*}
Here $H\left(  \tau\right)  ,\ \tau\in\left(  0,\infty\right)  ,\ $is given
by
\[
H\left(  \tau\right)  =\frac{1}{\left\vert B\left(  0,\tau\right)  \right\vert
}\int_{B\left(  0,\tau\right)  }\varphi\left(  y\right)  dy,\ \ \ \tau
\in\left(  0,\infty\right)
\]
and $B\left(  0,\tau\right)  $ is the open ball with radius $\tau>0$ centered
at the origin of $\mathbb{R}^{n}$.$\ $

\end{abstract}

\section{Introduction.\label{intro}  }

This article is a continuation of our previous ones \cite{TN, CT} and proves
several new interesting results.\ Consider the initial value problem%
\begin{equation}
\left\{
\begin{array}
[c]{l}%
u_{t}\left(  x,t\right)  =\bigtriangleup u\left(  x,t\right)  ,\ \ \ x\in
\mathbb{R}^{n},\ \ \ t>0,%
%TCIMACRO{\TeXButton{vspace}{\vspace{3mm}}}%
%BeginExpansion
\vspace{3mm}%
%EndExpansion
\\
u\left(  x,0\right)  =\varphi\left(  x\right)  ,\ \ \ x\in\mathbb{R}^{n},
\end{array}
\right.  \label{heat}%
\end{equation}
where $\varphi:\mathbb{R}^{n}\rightarrow\mathbb{R}$ is a given continuous
function. It is known that if $\varphi:\mathbb{R}^{n}\rightarrow\mathbb{R}$ is
a \textbf{continuous\ bounded} function, then the function given by
convolution integral
\begin{equation}
u\left(  x,t\right)  =\frac{1}{\left(  4\pi t\right)  ^{n/2}}\int
_{\mathbb{R}^{n}}e^{-\frac{\left\vert x-y\right\vert ^{2}}{4t}}\varphi\left(
y\right)  dy,\ \ \ x\in\mathbb{R}^{n},\ \ \ t>0 \label{CONV}%
\end{equation}
is a smooth solution of the heat equation on $\mathbb{R}^{n}\times
(0,\infty)\ $with $\lim_{\left(  x,t\right)  \rightarrow\left(  x_{0}%
,0\right)  }u\left(  x,t\right)  =\varphi\left(  x_{0}\right)  \ $for any
$x_{0}\in\mathbb{R}^{n}.\ $Due to the example by Tychonoff in 1935 (see the
book \cite{J}, Chapter 7), it is known that the Cauchy problem for the heat
equation (\ref{heat}) has no unique solution\ on $\mathbb{R}^{n}\times\left(
0,\infty\right)  \ $even if the initial data $\varphi\left(  x\right)  $ is
bounded (unless we impose certain growth condition of $u\left(  x,t\right)  $
as $\left\vert x\right\vert \rightarrow\infty$ for $t>0$).\ From now on, when
we say $u\left(  x,t\right)  $ is "the solution" of the heat equation with
$u\left(  x,0\right)  =\varphi\left(  x\right)  ,\ x\in\mathbb{R}^{n},$ we
always mean that it is the solution given by the convolution integral
(\ref{CONV}). As a consequence, if $\left\vert \varphi\right\vert \leq
M\ $on\ $\mathbb{R}^{n}\ $for\ some$\ $constant $M>0,\ $we also
have$\ \left\vert u\right\vert \leq M\ $on\ $\mathbb{R}^{n}\times
\lbrack0,\infty).$

In this paper, we shall always assume that $\varphi\in C^{0}\left(
\mathbb{R}^{n}\right)  \bigcap L^{\infty}\left(  \mathbb{R}^{n}\right)  \ $and
use$\ H\left(  \tau\right)  ,$ $\tau\in\lbrack0,\infty),$ to denote the
\textbf{average integral} of$\ \varphi\ $over the open ball $B\left(
0,\tau\right)  \subset\mathbb{R}^{n}$ centered at $x=0\ $with radius
$\tau>0,\ $i.e.%
\begin{equation}
H\left(  \tau\right)  =\frac{1}{\left\vert B\left(  0,\tau\right)  \right\vert
}\int_{B\left(  0,\tau\right)  }\varphi\left(  y\right)  dy=\frac{1}%
{\omega\left(  n\right)  \tau^{n}}\int_{B\left(  0,\tau\right)  }%
\varphi\left(  y\right)  dy,\ \ \ \tau\in\left(  0,\infty\right)  ,
\label{H-tau}%
\end{equation}
where $\omega\left(  n\right)  =\left\vert B\left(  0,1\right)  \right\vert $
is the volume of the unit ball in $\mathbb{R}^{n}.$ By continuity, if we
define $H\left(  0\right)  =\varphi\left(  0\right)  ,$ then$\ H\in
C^{0}[0,\infty)\bigcap L^{\infty}[0,\infty).$

The role played by $H\left(  \tau\right)  $ in the Cauchy problem (\ref{heat})
is the following beautiful \textbf{stabilization} result:

\begin{theorem}
\label{thm-K}(See \cite{E, K, RE}.) Assume $\varphi\in C^{0}\left(
\mathbb{R}^{n}\right)  \bigcap L^{\infty}\left(  \mathbb{R}^{n}\right)  $ and
$u\left(  x,t\right)  $ is the solution of the heat equation (\ref{heat})
given by (\ref{CONV}). Then
\begin{equation}
\lim_{t\rightarrow\infty}u\left(  0,t\right)  =0\ \ \ \text{if and only
if\ \ \ }\lim_{\tau\rightarrow\infty}\left(  \frac{1}{\left\vert B\left(
0,\tau\right)  \right\vert }\int_{B\left(  0,\tau\right)  }\varphi\left(
y\right)  dy\right)  =0. \label{K1}%
\end{equation}
Moreover,
\begin{equation}
\lim_{t\rightarrow\infty}\left(  \sup_{x\in\mathbb{R}^{n}}\left\vert u\left(
x,t\right)  \right\vert \right)  =0 \label{K2}%
\end{equation}
if and only if%
\begin{equation}
\lim_{\tau\rightarrow\infty}\left(  \sup_{x\in\mathbb{R}^{n}}\left\vert
\frac{1}{\left\vert B\left(  x,\tau\right)  \right\vert }\int_{B\left(
x,\tau\right)  }\varphi\left(  y\right)  dy\right\vert \right)  =0. \label{K3}%
\end{equation}
Here $\left\vert B\left(  x,\tau\right)  \right\vert $ is the volume of the
open ball in $\mathbb{R}^{n}$ centered at $x\in\mathbb{R}^{n}\ $with radius
$\tau>0.$
\end{theorem}

\begin{remark}
\label{rmk0}Similar result holds if we replace the limit value $0$ in
(\ref{K1}) by any real number $c.$
\end{remark}

\begin{remark}
\label{rmk-0}As demonstrated in \cite{CT}, Theorem \ref{thm-K} fails if the
initial data $\varphi\ $is \textbf{unbounded}. In this paper, we will always
assume $\varphi\in C^{0}\left(  \mathbb{R}^{n}\right)  \bigcap L^{\infty
}\left(  \mathbb{R}^{n}\right)  $.
\end{remark}

\begin{remark}
Since$\ \left\vert \varphi\right\vert \leq M\ $on\ $\mathbb{R}^{n}%
\ $for\ some$\ $constant $M>0,$ the convolution solution $u\left(  x,t\right)
$ satisfies the following gradient estimate
\begin{equation}
\left\vert \nabla_{x}u\left(  x,t\right)  \right\vert \leq\frac{c\left(
n\right)  M}{\sqrt{t}},\ \ \ \forall\ \left(  x,t\right)  \in\mathbb{R}%
^{n}\times\left(  0,\infty\right)  , \label{grad-est}%
\end{equation}
where $c\left(  n\right)  $ is a constant depending only on $n.\ $Hence if we
have$\ \lim_{t\rightarrow\infty}u\left(  0,t\right)  =0,$ we also have
$\lim_{t\rightarrow\infty}u\left(  x,t\right)  =0$ for all $x\in\mathbb{R}%
^{n}$ and the convergence is uniform in $x\in K$ where $K\subset\mathbb{R}%
^{n}$ is any compact set.
\end{remark}

By Theorem \ref{thm-K} and Remark \ref{rmk0}, the $\lim_{\tau\rightarrow
\infty}H\left(  \tau\right)  $ does not exist if and only if $\lim
_{t\rightarrow\infty}u\left(  0,t\right)  \ $does not exist. In this\ case, we
have the following comparison result concerning the oscillation of$\ u\left(
0,t\right)  $ and$\ H\left(  \tau\right)  $ as $t\rightarrow\infty
\ $and\ $\tau\rightarrow\infty\ $(see\ \cite{CT}), namely
\begin{equation}
\liminf_{\tau\rightarrow\infty}H\left(  \tau\right)  \leq\liminf
_{t\rightarrow\infty}u\left(  0,t\right)  <\limsup_{t\rightarrow\infty
}u\left(  0,t\right)  \leq\limsup_{\tau\rightarrow\infty}H\left(  \tau\right)
. \label{ineq-4}%
\end{equation}
The key point in the proof of (\ref{ineq-4}) is to use the
\textbf{representation formula}\ (see Lemma 4 in \cite{CT}):%
\begin{equation}
u\left(  0,t\right)  =\frac{2\omega\left(  n\right)  }{\pi^{n/2}}\int
_{0}^{\infty}e^{-z^{2}}z^{n+1}H\left(  \sqrt{4t}z\right)  dz,\ \ \ \forall
\ t\in\left(  0,\infty\right)  \label{rep}%
\end{equation}
together with the classical Fatou's Lemma.\ Same as in \cite{CT}, we shall use
$p\leq\alpha<\beta\leq q$ to denote the four finite limit values in
(\ref{ineq-4}) in this paper.

By analogy,\ can we also compare the oscillation of $u\left(  0,t\right)  ,$
$t\in\lbrack0,\infty),$ with the oscillation of the initial data
$\varphi\left(  x\right)  ?$ Since $\varphi\left(  x\right)  \ $is a function
defined on $\mathbb{R}^{n},$ it would be quite tricky unless we assume
$\varphi\left(  x\right)  =\varphi\left(  \tau\right)  ,\ \tau=\left\vert
x\right\vert \in\lbrack0,\infty),$ $x\in\mathbb{R}^{n},$ is a\textbf{ radial
function}. With this, both $u\left(  0,t\right)  $ and $\varphi\left(
\tau\right)  $ are one-variable functions defined on $[0,\infty)\ $and we can
compare their oscillations as $t\rightarrow\infty\ $and $\tau\rightarrow
\infty$ respectively. From now on, we shall assume that $\varphi\left(
x\right)  =\varphi\left(  \left\vert x\right\vert \right)  ,$ $x\in
\mathbb{R}^{n},$ is a\textbf{ }radial function.

Clearly, the limit$\ \lim_{\tau\rightarrow\infty}\varphi\left(  \tau\right)
=0\ $will imply $\lim_{t\rightarrow\infty}u\left(  0,t\right)  =0,$ but
\textbf{not} conversely.\ For example, the solution to the one-dimensional
heat equation $u_{t}=u_{xx}$ with radial initial data $\varphi\left(
x\right)  =\cos x,$ $x\in\lbrack0,\infty),$ is given by $u\left(  x,t\right)
=e^{-t}\cos x,$ which satisfies$\ \lim_{t\rightarrow\infty}u\left(
0,t\right)  =0.$ But we have $\liminf_{x\rightarrow\infty}\varphi\left(
x\right)  =-1,\ \limsup_{x\rightarrow\infty}\varphi\left(  x\right)  =1.\ $

For radial $\varphi\in C^{0}[0,\infty)\bigcap L^{\infty}[0,\infty),$ the
identities (\ref{CONV}) (at $x=0$)\ and (\ref{H-tau}) will become
\begin{equation}
u\left(  0,t\right)  =\frac{n\omega\left(  n\right)  }{\pi^{n/2}}\int
_{0}^{\infty}e^{-z^{2}}z^{n-1}\varphi\left(  \sqrt{4t}z\right)
dz,\ \ \ \forall\ t\in\left(  0,\infty\right)  \label{CONV-rad}%
\end{equation}
and
\begin{equation}
H\left(  \tau\right)  =\frac{n}{\tau^{n}}\int_{0}^{\tau}\varphi\left(
r\right)  r^{n-1}dr,\ \ \ \tau\in\left(  0,\infty\right)  ,\ \ \ H\left(
0\right)  =\varphi\left(  0\right)  . \label{H-tau-rad}%
\end{equation}
In particular, we note the similarity between (\ref{CONV-rad}) and (\ref{rep}).

Let $r=\liminf_{\tau\rightarrow\infty}\varphi\left(  \tau\right)
\ $and$\ s=\limsup_{\tau\rightarrow\infty}\varphi\left(  \tau\right)  .$ We
have the following general inequalities among the six numbers $\alpha
,\ \beta,\ p,\ q,\ r,\ s:$%
\begin{equation}
r\leq p\leq\alpha\leq\beta\leq q\leq s, \label{SIX}%
\end{equation}
where the first and last inequalities in (\ref{SIX}) can be derived from the
identity$\ H\left(  \tau\right)  =\frac{n}{\tau^{n}}\int_{0}^{\tau}%
\varphi\left(  r\right)  r^{n-1}dr$ and the definition of liminf and limsup.
Note that Theorem \ref{thm-K}\ says that $\alpha=\beta\ $if and only if $p=q.$
However, $\alpha=\beta$ does not necessarily imply $r=s\ $(as the above simple
example\ shows), i.e.\ the stabilization of $u\left(  0,t\right)  $ as
$t\rightarrow\infty$ does not imply the stabilization of $\varphi\left(
\tau\right)  $ as $\tau\rightarrow\infty.$

In this paper we are interested in prescribing the oscillation values
$\alpha<\beta\ $of $u\left(  0,t\right)  $ as $t\rightarrow\infty$ by choosing
suitable radial $\varphi\left(  \tau\right)  $ in (\ref{CONV-rad})\ or by
choosing suitable average integral $H\left(  \tau\right)  $ in (\ref{rep}) (by
(\ref{deri}) below, one can determine$\ \varphi\left(  \tau\right)  $ as long
as $H\left(  \tau\right)  $ is chosen). There are many choices for such
$\varphi\left(  \tau\right)  $ and $H\left(  \tau\right)  $.\ The interesting
results here are that one can also prescribe the oscillation values of
$\varphi\left(  \tau\right)  $ and $H\left(  \tau\right)  \ $as $\tau
\rightarrow\infty.$

More precisely, the results in this paper are:\ 

\begin{enumerate}
\item To prescribe, as $\tau\rightarrow\infty$ and $t\rightarrow\infty,$ the
oscillation behavior of $H\left(  \tau\right)  $ and $u\left(  0,t\right)
\ $for\textbf{ }four arbitrary finite numbers $p\leq\alpha<\beta\leq
q\ $satisfying the \textbf{symmetry condition}$\ $%
\begin{equation}
p+q=\alpha+\beta. \label{sym-cond}%
\end{equation}

\item To prescribe, as $\tau\rightarrow\infty$ and $t\rightarrow\infty,\ $the
oscillation behavior of $\varphi\left(  \tau\right)  $ and $u\left(
0,t\right)  $ for\ four\textbf{ }arbitrary finite numbers $r\leq\alpha
\leq\beta\leq s.$
\end{enumerate}

See Theorem \ref{thm-main-1}\ and Theorem \ref{thm-main-2} in the next section.

\section{Main results and their proofs.}

\subsection{Prescribing $H\left(  \tau\right)  $ and $u\left(  0,t\right)  .$}

The discussion of prescribing $H\left(  \tau\right)  $ and $u\left(
0,t\right)  \ $is motivated by Theorem \ref{thm-K}.\ Denote the four limit
values in (\ref{ineq-4}) as $p,\ \alpha,\ \beta,\ q$ respectively.\ In
\cite{CT},\ with the help of the formula (\ref{rep}), we have constructed two
bounded radial functions $\varphi\left(  x\right)  =\varphi\left(  \left\vert
x\right\vert \right)  ,$ $x\in\mathbb{R}^{n},$ so that the following are satisfied:

\begin{enumerate}
\item[(a).] $p=-1<\alpha=-\sqrt{A^{2}+B^{2}}<\beta=\sqrt{A^{2}+B^{2}}%
<q=1.\ $Here $A,\ B$ are the values of certain convergent improper integrals
satisfying\ $0<A^{2}+B^{2}<1$.\ 

\item[(b).] $p=\alpha<\beta=q\ $for arbitrary two finite numbers $\alpha
<\beta.$
\end{enumerate}

In the above two examples, we have the identity$\ $%
\begin{equation}
p+q=\alpha+\beta,\ \ \ \text{where}\ \ \ \alpha<\beta. \label{22}%
\end{equation}
Our new result in the following is that, as long as the symmetry condition
(\ref{22}) is satisfied, we can prescribe them. For four numbers $p\leq
\alpha<\beta\leq q$ to satisfy (\ref{22}), we must have either $p<\alpha
<\beta<q\ $or $p=\alpha<\beta=q.\ $Since the second case has been done, it
suffices to look at the first case.\ We have:

\begin{theorem}
\label{thm-main-1}Let $p<\alpha<\beta<q\ $be four \textbf{arbitrary} finite
numbers satisfying the identity $p+q=\alpha+\beta.\ $One can find a
\textbf{continuous bounded radial function}$\ \varphi\left(  x\right)
=\varphi\left(  \left\vert x\right\vert \right)  ,$ $x\in\mathbb{R}^{n},\ $so
that its average integral $H\left(  \tau\right)  \ $and the convolution
solution (\ref{CONV}) of the problem (\ref{heat})\ satisfy
\begin{equation}
\liminf_{\tau\rightarrow\infty}H\left(  \tau\right)  =p<\liminf_{t\rightarrow
\infty}u\left(  0,t\right)  =\alpha<\limsup_{t\rightarrow\infty}u\left(
0,t\right)  =\beta<\limsup_{\tau\rightarrow\infty}H\left(  \tau\right)
=q.\label{four}%
\end{equation}
More precisely, we can choose the initial data $\varphi\left(  x\right)
=\varphi\left(  \left\vert x\right\vert \right)  ,$ $x\in\mathbb{R}^{n},$ as%
\begin{equation}
\varphi\left(  \tau\right)  =\frac{q-p}{2}\sin\left(  m\log\left(
\tau+1\right)  \right)  +\frac{q-p}{2}\frac{m\tau}{n\left(  \tau+1\right)
}\cos\left(  m\log\left(  \tau+1\right)  \right)  +\frac{q+p}{2},\label{phi}%
\end{equation}
where$\ \tau\in\lbrack0,\infty)$ and$\ m\in\left(  0,\infty\right)  $ is a
number satisfying the identity%
\begin{align}
&  \left(  \frac{2\omega\left(  n\right)  }{\pi^{n/2}}\int_{0}^{\infty
}e^{-z^{2}}z^{n+1}\cos\left(  m\log z\right)  dz\right)  ^{2}+\left(
\frac{2\omega\left(  n\right)  }{\pi^{n/2}}\int_{0}^{\infty}e^{-z^{2}}%
z^{n+1}\sin\left(  m\log z\right)  dz\right)  ^{2}\nonumber\\
&  =\left(  \frac{\beta-\alpha}{q-p}\right)  ^{2}\in\left(  0,1\right)
.\label{IVT}%
\end{align}

\end{theorem}

\begin{remark}
Without the condition $p+q=\alpha+\beta,$ we \textbf{do not know} how to
achieve the result (\ref{four}) in general. However, we can construct some
specific example with $p+q\neq\alpha+\beta.$ See Section \ref{sec-exam}.
\end{remark}

\begin{remark}
\label{rmk-rad}Note that if the initial data $\varphi\left(  x\right)
=\varphi\left(  \left\vert x\right\vert \right)  $ is radial, the convolution
solution $u\left(  x,t\right)  \ $of the heat equation will also be radial in
$x\in\mathbb{R}^{n}$ for each fixed time $t>0$.\ 
\end{remark}

\subsection{Prescribing $\varphi\left(  \tau\right)  $ and $u\left(
0,t\right)  .$}

It is due to Theorem \ref{thm-K} that we study the oscillation behavior
of$\ H\left(  \tau\right)  $ and its effect on $u\left(  0,t\right)  .$
Instead, since $\varphi\left(  \tau\right)  $ is the initial data, it is
perhaps more natural to study the oscillation relation between$\ \varphi
\left(  \tau\right)  $ and $u\left(  0,t\right)  .$ However, there is a major
difference here.\ Since we always assume that $\varphi\left(  \tau\right)  \in
C^{0}[0,\infty)\bigcap L^{\infty}[0,\infty)$ is continuous and bounded, its
average integral function $H\left(  \tau\right)  $ will also be continuous and
bounded and satisfy $H^{\prime}\left(  \tau\right)  =O\left(  1/\tau\right)  $
as $\tau\rightarrow\infty,$ which we call it a \textbf{slow oscillation} if we
have $p<q\ $(see (\ref{d-est}) and Remark \ref{rmk1} below).\ Any such
oscillation will cause $u\left(  0,t\right)  $ to have a slow oscillation too
(by Theorem \ref{thm-K} we have $p<q\ $if and only if $\alpha<\beta$).

On the other hand, $\varphi\left(  \tau\right)  \in C^{0}[0,\infty)\bigcap
L^{\infty}[0,\infty)$ may not have a slow oscillation (i.e. it may not satisfy
$\varphi^{\prime}\left(  \tau\right)  =O\left(  1/\tau\right)  \ $as
$\tau\rightarrow\infty\ $even if it is differentiable) and so its oscillation
(with $r<s$)\ may not pass to $u\left(  0,t\right)  $.\ In fact, for any
$2\pi$\textbf{-periodic} regular oscillation function $\varphi\left(
\tau\right)  ,$ $\tau\in\lbrack0,\infty),$ the Cauchy problem (\ref{heat})
with$\ u\left(  x,0\right)  =\varphi\left(  \left\vert x\right\vert \right)
,$ $x\in\mathbb{R}^{n},\ $will have\ the convergence$\ \lim_{t\rightarrow
\infty}u\left(  0,t\right)  =\frac{1}{2\pi}\int_{0}^{2\pi}\varphi\left(
\tau\right)  d\tau$ (see Lemma \ref{lem-ave} below). However, due to the
similarity between the two representation formulas (\ref{rep}) and
(\ref{CONV-rad}), if $\varphi\left(  \tau\right)  \ $has a slow
oscillation\textbf{ }as $\tau\rightarrow\infty\ $with $r<s,$ then it has the
chance of passing to the function $u\left(  0,t\right)  \ $as $t\rightarrow
\infty$ (see Lemma \ref{lem-slow} below).

The major result for prescribing $\varphi\left(  \tau\right)  $ and $u\left(
0,t\right)  \ $is the following:

\begin{theorem}
\label{thm-main-2}Let $r<\alpha<\beta<s\ $be four \textbf{arbitrary} finite
numbers. One can find a \textbf{continuous bounded radial function }%
$\varphi\left(  x\right)  =\varphi\left(  \left\vert x\right\vert \right)  ,$
$x\in\mathbb{R}^{n},$ so that it and the convolution solution (\ref{CONV}) of
the problem (\ref{heat})\ satisfy%
\begin{equation}
\liminf_{\tau\rightarrow\infty}\varphi\left(  \tau\right)  =r<\liminf
_{t\rightarrow\infty}u\left(  0,t\right)  =\alpha<\limsup_{t\rightarrow\infty
}u\left(  0,t\right)  =\beta<\limsup_{\tau\rightarrow\infty}\varphi\left(
\tau\right)  =s.
\end{equation}
More precisely, we can choose the initial data $\varphi\left(  x\right)
=\varphi\left(  \left\vert x\right\vert \right)  ,$ $x\in\mathbb{R}^{n},$ as%
\begin{equation}
\varphi\left(  \tau\right)  =\left[  C_{1}\sin\left(  m\log\left(
\tau+1\right)  \right)  +C_{2}\right]  +\psi\left(  \tau\right)
,\ \ \ \tau\in\lbrack0,\infty), \label{phi-0}%
\end{equation}
for some suitable constants $m>0,\ C_{1}>0,\ C_{2}\ $depending on
$r,\ \alpha,\ \beta,\ s\ $and some $2\pi$-periodic continuous
nonzero\ function $\psi\left(  \tau\right)  $ satisfying$\ \frac{1}{2\pi}%
\int_{0}^{2\pi}\psi\left(  \tau\right)  d\tau=0.$
\end{theorem}

\begin{remark}
Theorem \ref{thm-main-2} is actually valid for four arbitrary finite
numbers$\ r\leq\alpha\leq\beta\leq s.\ $We will deal with the remaining cases
in Theorem \ref{thm-main-3}\ below.
\end{remark}

\begin{remark}
This is a comparison between\ Theorem \ref{thm-main-1}\ and Theorem
\ref{thm-main-2}.\ For any initial data $\varphi\in C^{0}\left(
\mathbb{R}^{n}\right)  \bigcap L^{\infty}\left(  \mathbb{R}^{n}\right)  ,$ no
matter it is radial or not, the function $H\left(  \tau\right)  $ is always a
slow oscillation function as $\tau\rightarrow\infty$. Therefore, the function
space for $H\left(  \tau\right)  $ is smaller than the function space for
$\varphi\left(  \tau\right)  ,\ $which may explain why we have more general
result in Theorem \ref{thm-main-2}. Note that the function $\varphi\left(
\tau\right)  $ in (\ref{phi-0}) does not have a\ slow oscillation as
$\tau\rightarrow\infty$ due to the term $\psi\left(  \tau\right)  $.
\end{remark}

\subsubsection{Proof\ of Theorem \ref{thm-main-1}.}

The proof is an interesting modification of the argument in Section 3.1.1 of
\cite{CT}, which can achieve much more general result.\ Some motivation should
be explained at the beginning.\ First, we choose $\varphi\left(  x\right)
=\varphi\left(  \left\vert x\right\vert \right)  $ to be a radial function in
order to simplify the computation. We think a non-radial function can also be
found as long as we can overcome the computational complexity. Second, we
choose the average integral function $H\left(  \tau\right)  $ first and then
go back to find its corresponding $\varphi\left(  \left\vert x\right\vert
\right)  .\ $The function $H\left(  \tau\right)  \ $has to be
\textbf{oscillatory} so that we have$\ \liminf_{\tau\rightarrow\infty}H\left(
\tau\right)  <\limsup_{\tau\rightarrow\infty}H\left(  \tau\right)  .$
Moreover, since we require $\varphi$ to bounded and continuous$\ $(in view of
Theorem \ref{thm-K}),\ its average integral $H\left(  \tau\right)  $ will
satisfy the\ derivative estimate$\ H^{\prime}\left(  \tau\right)  =O\left(
1/\tau\right)  \ $as$\ \tau\rightarrow\infty\ $(see (\ref{d-est})
below).\ Therefore, a natural choice is roughly like the function $\sin\left(
\log\tau\right)  ,\ \tau\in\left(  0,\infty\right)  .\ $After suitable
modification of the function $\sin\left(  \log\tau\right)  ,$ the proof can be achieved.

Given $\varphi\left(  x\right)  =\varphi\left(  \left\vert x\right\vert
\right)  \in C^{0}[0,\infty)\bigcap L^{\infty}[0,\infty),$ its average
integral function
\begin{equation}
H\left(  \tau\right)  =\frac{1}{\omega\left(  n\right)  \tau^{n}}%
\int_{B\left(  0,\tau\right)  }\varphi\left(  y\right)  dy=\frac{n}{\tau^{n}%
}\int_{0}^{\tau}\varphi\left(  r\right)  r^{n-1}dr,\ \ \ \tau\in\left(
0,\infty\right)  ,\ \ \ H\left(  0\right)  =\varphi\left(  0\right)  ,
\end{equation}
is a bounded continuous function on$\ [0,\infty),$ differentiable in $\left(
0,\infty\right)  ,$ and satisfies the identity
\begin{equation}
H^{\prime}\left(  \tau\right)  =-\frac{n}{\tau}H\left(  \tau\right)  +\frac
{n}{\tau}\varphi\left(  \tau\right)  ,\ \ \ \forall\ \tau\in\left(
0,\infty\right)  . \label{deri}%
\end{equation}
Hence it has the asymptotic behavior $\ $%
\begin{equation}
H^{\prime}\left(  \tau\right)  =O\left(  \frac{1}{\tau}\right)
\ \ \ \text{as\ \ \ }\tau\rightarrow\infty. \label{d-est}%
\end{equation}
Moreover, the radial function $\varphi\left(  \tau\right)  $ can be determined
from$\ H\left(  \tau\right)  \ $by the identity (\ref{deri}).

Let $m\in\left(  0,\infty\right)  \ $be a fixed number. Motivated by our
principle stated above, we temporarily choose $H\left(  \tau\right)  $ (will
modify it later on)\ to be equal to
\begin{equation}
H\left(  \tau\right)  =\frac{n}{\tau^{n}}\int_{0}^{\tau}\varphi\left(
r\right)  r^{n-1}dr=\sin\left(  m\log\left(  \tau+1\right)  \right)
,\ \ \ \tau\in\left(  0,\infty\right)  ,\ \ \ H\left(  0\right)  =0,
\label{Hm-tau}%
\end{equation}
which, by (\ref{deri}), gives
\begin{equation}
\varphi\left(  \tau\right)  =\sin\left(  m\log\left(  \tau+1\right)  \right)
+\frac{m\tau}{n\left(  \tau+1\right)  }\cos\left(  m\log\left(  \tau+1\right)
\right)  ,\ \ \ \tau\in\left(  0,\infty\right)  . \label{3-tau}%
\end{equation}
That is, if we choose $\varphi\left(  \tau\right)  $ to be\ the function given
by (\ref{3-tau}) and define $\varphi\left(  0\right)  =0,\ $then the function
$\varphi\left(  x\right)  =\varphi\left(  \left\vert x\right\vert \right)  ,$
$x\in\mathbb{R}^{n},\ $will be a radial function defined on $\mathbb{R}%
^{n},\ $lying in the space $C^{0}\left(  \mathbb{R}^{n}\right)  \bigcap
L^{\infty}\left(  \mathbb{R}^{n}\right)  ,\ $and its average integral
is$\ $given\ by$\ H\left(  \tau\right)  =\sin\left(  m\log\left(
\tau+1\right)  \right)  $ for all $\tau\in\lbrack0,\infty).\ $In particular,
we have $p=-1,\ q=1.$

To find the values of $\alpha,\ \beta,\ $we use the representation formula
(\ref{rep}).\ We have
\begin{align}
&  u\left(  0,t\right)  =\frac{2\omega\left(  n\right)  }{\pi^{n/2}}\int
_{0}^{\infty}e^{-z^{2}}z^{n+1}\sin\left(  m\log\left(  \sqrt{4t}z+1\right)
\right)  dz\nonumber\\
&  =\left\{
\begin{array}
[c]{l}%
\left[  \frac{2\omega\left(  n\right)  }{\pi^{n/2}}\int_{0}^{\infty}e^{-z^{2}%
}z^{n+1}\cos\left(  m\log\left(  z+\frac{1}{\sqrt{4t}}\right)  \right)
dz\right]  \sin\left(  m\log\sqrt{4t}\right)
%TCIMACRO{\TeXButton{vspace}{\vspace{3mm}}}%
%BeginExpansion
\vspace{3mm}%
%EndExpansion
\\
+\left[  \frac{2\omega\left(  n\right)  }{\pi^{n/2}}\int_{0}^{\infty}%
e^{-z^{2}}z^{n+1}\sin\left(  m\log\left(  z+\frac{1}{\sqrt{4t}}\right)
\right)  dz\right]  \cos\left(  m\log\sqrt{4t}\right)  ,\ \ \ t\in\left(
0,\infty\right)
\end{array}
\right.  \ \label{u0t-for}%
\end{align}
and the Lebesgue Dominated Convergence Theorem implies%
\begin{align}
&  \lim_{t\rightarrow\infty}\frac{2\omega\left(  n\right)  }{\pi^{n/2}}%
\int_{0}^{\infty}e^{-z^{2}}z^{n+1}\cos\left(  m\log\left(  z+\frac{1}%
{\sqrt{4t}}\right)  \right)  dz\nonumber\\
&  =\frac{2\omega\left(  n\right)  }{\pi^{n/2}}\int_{0}^{\infty}e^{-z^{2}%
}z^{n+1}\cos\left(  m\log z\right)  dz:=A\left(  m\right)  \in\left(
-1,1\right)  \label{Am}%
\end{align}
and
\begin{align}
&  \lim_{t\rightarrow\infty}\frac{2\omega\left(  n\right)  }{\pi^{n/2}}%
\int_{0}^{\infty}e^{-z^{2}}z^{n+1}\sin\left(  m\log\left(  z+\frac{1}%
{\sqrt{4t}}\right)  \right)  dz\nonumber\\
&  =\frac{2\omega\left(  n\right)  }{\pi^{n/2}}\int_{0}^{\infty}e^{-z^{2}%
}z^{n+1}\sin\left(  m\log z\right)  dz:=B\left(  m\right)  \in\left(
-1,1\right)  , \label{Bm}%
\end{align}
where we can use the identity$\ \frac{2\omega\left(  n\right)  }{\pi^{n/2}%
}\int_{0}^{\infty}e^{-z^{2}}z^{n+1}dz=1\ $to know that $A\left(  m\right)
,\ B\left(  m\right)  \in\left(  -1,1\right)  .\ $By H\"{o}lder inequality,
the constant $A\left(  m\right)  $ satisfies the estimate
\begin{align*}
A^{2}\left(  m\right)   &  =\left(  \frac{2\omega\left(  n\right)  }{\pi
^{n/2}}\int_{0}^{\infty}\sqrt{e^{-z^{2}}z^{n+1}}\sqrt{e^{-z^{2}}z^{n+1}}%
\cos\left(  m\log z\right)  dz\right)  ^{2}\\
&  <\left(  \frac{2\omega\left(  n\right)  }{\pi^{n/2}}\int_{0}^{\infty
}e^{-z^{2}}z^{n+1}dz\right)  \left(  \frac{2\omega\left(  n\right)  }%
{\pi^{n/2}}\int_{0}^{\infty}e^{-z^{2}}z^{n+1}\cos^{2}\left(  m\log z\right)
dz\right) \\
&  =\frac{2\omega\left(  n\right)  }{\pi^{n/2}}\int_{0}^{\infty}e^{-z^{2}%
}z^{n+1}\cos^{2}\left(  m\log z\right)  dz
\end{align*}
and similarly
\[
B^{2}\left(  m\right)  <\frac{2\omega\left(  n\right)  }{\pi^{n/2}}\int
_{0}^{\infty}e^{-z^{2}}z^{n+1}\sin^{2}\left(  m\log z\right)  dz
\]
and we conclude%
\begin{equation}
0<A^{2}\left(  m\right)  +B^{2}\left(  m\right)  <\frac{2\omega\left(
n\right)  }{\pi^{n/2}}\int_{0}^{\infty}e^{-z^{2}}z^{n+1}\left[  \cos
^{2}\left(  m\log z\right)  +\sin^{2}\left(  m\log z\right)  \right]  dz=1.
\end{equation}
Therefore, we have the convergence
\begin{equation}
\lim_{t\rightarrow\infty}\left\vert u\left(  0,t\right)  -\left[  A\left(
m\right)  \sin\left(  m\log\sqrt{4t}\right)  +B\left(  m\right)  \cos\left(
m\log\sqrt{4t}\right)  \right]  \right\vert =0, \label{u0t-limit}%
\end{equation}
which gives the conclusion%
\begin{align}
-1  &  =\liminf_{\tau\rightarrow\infty}H\left(  \tau\right)  <-\sqrt
{A^{2}\left(  m\right)  +B^{2}\left(  m\right)  }=\liminf_{t\rightarrow\infty
}u\left(  0,t\right) \nonumber\\
&  <\ \limsup_{t\rightarrow\infty}u\left(  0,t\right)  =\sqrt{A^{2}\left(
m\right)  +B^{2}\left(  m\right)  }<\limsup_{\tau\rightarrow\infty}H\left(
\tau\right)  =1 \label{1AmBm1}%
\end{align}
for fixed $m\in\left(  0,\infty\right)  .$

The idea next is to view$\ m\in\left(  0,\infty\right)  $ as a variable.\ We
observe the following:

\begin{lemma}
\label{lem-m}We have
\begin{equation}
\lim_{m\rightarrow\infty}A\left(  m\right)  =\lim_{m\rightarrow\infty}B\left(
m\right)  =0,\ \ \ m\in\left(  0,\infty\right)  . \label{AmBm}%
\end{equation}

\end{lemma}

%

%TCIMACRO{\TeXButton{Proof}{\proof}}%
%BeginExpansion
\proof
%EndExpansion
For any $\varepsilon>0,\ $one can find small number$\ \delta>0\ $and large
number$\ M>0,\ $both are independent\ of $m\in\left(  0,\infty\right)  ,$ such
that%
\begin{equation}
\left\vert \int_{0}^{\delta}e^{-z^{2}}z^{n+1}\cos\left(  m\log z\right)
dz\right\vert +\left\vert \int_{M}^{\infty}e^{-z^{2}}z^{n+1}\cos\left(  m\log
z\right)  dz\right\vert <\varepsilon. \label{Am-1}%
\end{equation}
On the other hand, by the change of variables, we have
\[
\int_{\delta}^{M}e^{-z^{2}}z^{n+1}\cos\left(  m\log z\right)  dz=\int
_{\log\delta}^{\log M}F\left(  x\right)  \cos\left(  mx\right)
dx,\ \ \ x=\log z,
\]
where$\ F\left(  x\right)  =e^{-e^{2x}+\left(  n+2\right)  x}.\ $The familiar
\textbf{Riemann-Lebesgue lemma }in analysis (or just use integration by
parts)\ implies
\begin{equation}
\lim_{m\rightarrow\infty}\int_{\log\delta}^{\log M}F\left(  x\right)
\cos\left(  mx\right)  dx=0. \label{Am-2}%
\end{equation}
which, together with (\ref{Am-1}), implies$\ \lim_{m\rightarrow\infty}A\left(
m\right)  =0.\ $The proof of $\lim_{m\rightarrow\infty}B\left(  m\right)
=0\ $is similar.$%
%TCIMACRO{\TeXButton{hfill}{\hfill}}%
%BeginExpansion
\hfill
%EndExpansion
\square$

\ \ \ \ 

Since estimate\ (\ref{Am-1}) is valid for all $m\in\left(  0,\infty\right)
,\ $the two improper integrals%
\begin{equation}
\int_{0}^{\infty}e^{-z^{2}}z^{n+1}\cos\left(  m\log z\right)  dz,\ \ \ \int
_{0}^{\infty}e^{-z^{2}}z^{n+1}\sin\left(  m\log z\right)  dz \label{2-int}%
\end{equation}
\textbf{converge uniformly} with respect to $m\in\left(  0,\infty\right)
.\ $As a consequence, both $A\left(  m\right)  \ $and $B\left(  m\right)
\ $are \textbf{continuous functions }of $m\in\left(  0,\infty\right)  ,\ $with%
\begin{equation}
\lim_{m\rightarrow\infty}A\left(  m\right)  =\lim_{m\rightarrow\infty}B\left(
m\right)  =0,\ \ \ \lim_{m\rightarrow0}A\left(  m\right)  =1,\ \ \ \lim
_{m\rightarrow0}B\left(  m\right)  =0, \label{AmBm-1}%
\end{equation}
where the last two limits in (\ref{AmBm-1}) are due to the Lebesgue Dominated
Convergence Theorem.$\allowbreak$ Moreover, the intermediate value theorem
implies the existence of a number $m\in\left(  0,\infty\right)  \ $satisfying
\begin{equation}
\sqrt{A^{2}\left(  m\right)  +B^{2}\left(  m\right)  }=\frac{\beta-\alpha
}{q-p}\in\left(  0,1\right)  . \label{ivt}%
\end{equation}

Now we choose $H\left(  \tau\right)  $ as
\begin{equation}
H\left(  \tau\right)  =\frac{q-p}{2}\sin\left(  m\log\left(  \tau+1\right)
\right)  +\frac{q+p}{2},\ \ \ H\left(  0\right)  =\frac{q+p}{2},\ \ \ \tau
\in\lbrack0,\infty), \label{H-tau-m}%
\end{equation}
where $m\in\left(  0,\infty\right)  \ $is the number satisfying (\ref{ivt}%
).\ By (\ref{Hm-tau}) and (\ref{3-tau}), its corresponding initial bounded
radial function $\varphi\left(  x\right)  =\varphi\left(  \left\vert
x\right\vert \right)  ,$ $x\in\mathbb{R}^{n},$ is given by%
\begin{equation}
\varphi\left(  \tau\right)  =\left\{
\begin{array}
[c]{l}%
\frac{q-p}{2}\sin\left(  m\log\left(  \tau+1\right)  \right)
%TCIMACRO{\TeXButton{vspace}{\vspace{3mm}}}%
%BeginExpansion
\vspace{3mm}%
%EndExpansion
\\
+\frac{q-p}{2}\frac{m\tau}{n\left(  \tau+1\right)  }\cos\left(  m\log\left(
\tau+1\right)  \right)  +\frac{q+p}{2},\ \ \ \varphi\left(  0\right)
=\frac{q+p}{2},\ \ \ \tau\in\lbrack0,\infty).
\end{array}
\right.  \label{phi-1}%
\end{equation}
\bigskip Similar to (\ref{u0t-limit}), we have
\begin{equation}
\lim_{t\rightarrow\infty}\left\vert u\left(  0,t\right)  -\left\{  \frac
{q-p}{2}\left[  A\left(  m\right)  \sin\left(  m\log\sqrt{4t}\right)
+B\left(  m\right)  \cos\left(  m\log\sqrt{4t}\right)  \right]  +\frac{q+p}%
{2}\right\}  \right\vert =0,
\end{equation}
which gives%
\[
\liminf_{t\rightarrow\infty}u\left(  0,t\right)  =-\frac{q-p}{2}\sqrt
{A^{2}\left(  m\right)  +B^{2}\left(  m\right)  }+\frac{q+p}{2}=-\frac
{\beta-\alpha}{2}+\frac{\beta+\alpha}{2}=\alpha
\]
and%
\[
\limsup_{t\rightarrow\infty}u\left(  0,t\right)  =\frac{q-p}{2}\sqrt
{A^{2}\left(  m\right)  +B^{2}\left(  m\right)  }+\frac{q+p}{2}=\frac
{\beta-\alpha}{2}+\frac{\beta+\alpha}{2}=\beta.
\]
Since we clearly have$\ \liminf_{\tau\rightarrow\infty}H\left(  \tau\right)
=p,\ \limsup_{\tau\rightarrow\infty}H\left(  \tau\right)  =q,$ the proof of
Theorem \ref{thm-main-1} is now complete.$%
%TCIMACRO{\TeXButton{hfill}{\hfill}}%
%BeginExpansion
\hfill
%EndExpansion
\square${}

\begin{remark}
The proof\ of Theorem \ref{thm-main-1} also reveals the following interesting
observation. For fixed $q-p>0,\ $if we have small $\beta-\alpha,\ $then by
(\ref{ivt}) and (\ref{AmBm}), we will have large$\ m\in\left(  0,\infty
\right)  \ $and the function $\varphi\left(  \tau\right)  $ in (\ref{phi-1})
will tend to be unbounded due to the term$\ m\tau/\left(  n\left(
\tau+1\right)  \right)  .$ This may suggest that we cannot find a bounded
radial function$\ \varphi\left(  \tau\right)  $ satisfying $p<\alpha=\beta<q.$
This matches with the result in Theorem \ref{thm-K}, which implies that if we
have$\ p<\alpha=\beta<q,$ then the initial data $\varphi\in C^{0}\left(
\mathbb{R}^{n}\right)  \ $must be unbounded. \ 
\end{remark}

\begin{remark}
By (\ref{grad-est}), for a given initial data $\varphi\in C^{0}\left(
\mathbb{R}^{n}\right)  \bigcap L^{\infty}\left(  \mathbb{R}^{n}\right)  ,$ if
we have$\ $%
\[
\liminf_{t\rightarrow\infty}u\left(  0,t\right)  =\alpha<\limsup
_{t\rightarrow\infty}u\left(  0,t\right)  =\beta,
\]
then for any fixed $x\in\mathbb{R}^{n}\ $we also have%
\begin{equation}
\liminf_{t\rightarrow\infty}u\left(  x,t\right)  =\alpha<\limsup
_{t\rightarrow\infty}u\left(  x,t\right)  =\beta. \label{too}%
\end{equation}

\end{remark}

\begin{remark}
\label{rmk1}It has been shown in \cite{CT} (see its equations (16)\ and
(37))\ that both functions $H\left(  \tau\right)  \ $and$\ u\left(
0,t\right)  $ satisfy%
\begin{equation}
H^{\prime}\left(  \tau\right)  =O\left(  \frac{1}{\tau}\right)
\ \ \text{as\ \ }\tau\rightarrow\infty,\ \ \ \ \ u_{t}\left(  0,t\right)
=O\left(  \frac{1}{t}\right)  \ \ \text{as\ \ }t\rightarrow\infty, \label{OO}%
\end{equation}
as long as the initial data $\varphi\ $lies in the space$\ C^{0}\left(
\mathbb{R}^{n}\right)  \bigcap L^{\infty}\left(  \mathbb{R}^{n}\right)  $ (no
matter it is radial or not). By (\ref{OO}), we may say that any oscillation
in$\ H\left(  \tau\right)  \ $(with $p<q$)\ or in $u\left(  0,t\right)  $
(with $\alpha<\beta$)\ is a\textbf{ slow oscillation. \ }
\end{remark}

As a consequence of Theorem \ref{thm-main-1}, we have the following two
corollaries, which say that we can \textbf{prescribe} the oscillation of
$H\left(  \tau\right)  \ $and $u\left(  0,t\right)  $ for \textbf{three}
arbitrary different numbers.\ 

\begin{corollary}
\label{cor-2}For any three different numbers, denoted as $p<\alpha<\beta
,\ $one can find a continuous\ bounded\ radial function$\ \varphi\left(
x\right)  =\varphi\left(  \left\vert x\right\vert \right)  ,$ $x\in
\mathbb{R}^{n},\ $so that its average integral and the solution (\ref{CONV}%
)\ of the problem\ (\ref{heat}) satisfies
\begin{equation}
\liminf_{\tau\rightarrow\infty}H\left(  \tau\right)  =p<\liminf_{t\rightarrow
\infty}u\left(  0,t\right)  =\alpha<\limsup_{t\rightarrow\infty}u\left(
0,t\right)  =\beta<\limsup_{\tau\rightarrow\infty}H\left(  \tau\right)  =q,
\label{pabq}%
\end{equation}
where $q=\alpha+\beta-p.\ $Similarly, for any three different\ numbers,
denoted as $\alpha<\beta<q,\ $one can find a bounded\ radial function\textbf{
}satisfying (\ref{pabq}), where now $p=\alpha+\beta-q.\ $
\end{corollary}

\begin{corollary}
For any\textbf{ }three\textbf{ }different\ numbers, denoted as $p<\alpha
<q,\ $with $\alpha<\left(  p+q\right)  /2,\ $the same result as in Corollary
\ref{cor-2} holds, where now $\beta=p+q-\alpha.\ $Similarly, for any three
different numbers, denoted as $p<\beta<q,\ $with $\beta>\left(  p+q\right)
/2,\ $the same result as in Corollary \ref{cor-2} holds, where now
$\alpha=p+q-\beta.\ $
\end{corollary}

\subsubsection{Proof\ of Theorem \ref{thm-main-2}.}

To prove Theorem \ref{thm-main-2}, we first need the following:

\begin{lemma}
\label{lem-ave}Assume $\varphi\left(  \tau\right)  \ $is a $2\pi
$\textbf{-periodic} radial function defined on $\tau\in\lbrack0,\infty)$%
.$\ $We have the convergence
\begin{equation}
\lim_{t\rightarrow\infty}u\left(  0,t\right)  =\lim_{\tau\rightarrow\infty
}H\left(  \tau\right)  =\lim_{\tau\rightarrow\infty}\frac{n}{\tau^{n}}\int
_{0}^{\tau}\varphi\left(  r\right)  r^{n-1}dr=\frac{1}{2\pi}\int_{0}^{2\pi
}\varphi\left(  \tau\right)  d\tau. \label{ave-1}%
\end{equation}

\end{lemma}%

%TCIMACRO{\TeXButton{Proof}{\proof}}%
%BeginExpansion
\proof
%EndExpansion
Since $\varphi\in C^{0}[0,\infty)\bigcap L^{\infty}[0,\infty),$ by Remark
\ref{rmk0} in Theorem \ref{thm-K}, it suffices to prove the identity for
$H\left(  \tau\right)  .\ $The proof is quite straightforward for the case
$n=1\ $(see Remark 6 in \cite{TN}), but may need a trick for $n>1.\ $For large
$\tau>0,$ we can express it as $\tau=2m\pi+R\ $for some $m\in\mathbb{N}$ and
$R\in\lbrack0,2\pi),$ with the understanding that both $m\ $and $R\ $depend on
$\tau$ and as $\tau\rightarrow\infty$ we have $m\rightarrow\infty,$ hence we
obtain
\begin{align}
&  \frac{n}{\tau^{n}}\int_{0}^{\tau}\varphi\left(  r\right)  r^{n-1}%
dr\nonumber\\
&  =\frac{n}{\left(  2m\pi+R\right)  ^{n}}\left(
\begin{array}
[c]{l}%
\int_{0}^{2\pi}\varphi\left(  r\right)  r^{n-1}dr+\int_{2\pi}^{4\pi}%
\varphi\left(  r\right)  r^{n-1}dr+\int_{4\pi}^{6\pi}\varphi\left(  r\right)
r^{n-1}dr\\
+\cdot\cdot\cdot+\int_{2\left(  m-1\right)  \pi}^{2m\pi}\varphi\left(
r\right)  r^{n-1}dr+\int_{2m\pi}^{2m\pi+R}\varphi\left(  r\right)  r^{n-1}dr
\end{array}
\right)  . \label{sum}%
\end{align}
If we do the change of variables%
\[
\int_{2\pi}^{4\pi}\varphi\left(  r\right)  r^{n-1}dr=\int_{0}^{2\pi}%
\varphi\left(  s\right)  \left(  s+2\pi\right)  ^{n-1}ds,\ \ \ r=s+2\pi,
\]
and%
\[
\int_{4\pi}^{6\pi}\varphi\left(  r\right)  r^{n-1}dr=\int_{0}^{2\pi}%
\varphi\left(  s\right)  \left(  s+4\pi\right)  ^{n-1}ds,\ \ \ r=s+4\pi,
\]
$\ ...,\ $etc.,\ (\ref{sum}) becomes
\begin{align}
&  \frac{n}{\tau^{n}}\int_{0}^{\tau}\varphi\left(  r\right)  r^{n-1}%
dr\nonumber\\
&  =\frac{n}{\left(  2m\pi+R\right)  ^{n}}\left(
\begin{array}
[c]{l}%
\int_{0}^{2\pi}\varphi\left(  s\right)  \left[  s^{n-1}+\left(  s+2\pi\right)
^{n-1}+\left(  s+4\pi\right)  ^{n-1}+\cdot\cdot\cdot+\left(  s+2\left(
m-1\right)  \pi\right)  ^{n-1}\right]  ds%
%TCIMACRO{\TeXButton{vspace}{\vspace{3mm}}}%
%BeginExpansion
\vspace{3mm}%
%EndExpansion
\\
+\int_{0}^{R}\varphi\left(  s\right)  \left(  s+2m\pi\right)  ^{n-1}ds
\end{array}
\right) \nonumber\\
&  :=I+II. \label{ave-2}%
\end{align}
For the second term $II\ $in (\ref{ave-2}), we have%
\begin{equation}
\lim_{m\rightarrow\infty}\left\vert II\right\vert \leq\lim_{m\rightarrow
\infty}\left(  \frac{n}{\left(  2m\pi+R\right)  ^{n}}\cdot2\pi\max
_{s\in\left[  0,2\pi\right]  }\left\vert \varphi\left(  s\right)  \right\vert
\cdot\left(  2\pi+2m\pi\right)  ^{n-1}\right)  =0. \label{ave-3}%
\end{equation}
To estimate the first term $I\ $in (\ref{ave-2}), we first look at the
integral
\begin{align*}
&  \frac{n}{\left(  2m\pi\right)  ^{n}}\int_{0}^{2\pi}\varphi\left(  s\right)
\left[  s^{n-1}+\left(  s+2\pi\right)  ^{n-1}+\left(  s+4\pi\right)
^{n-1}+\cdot\cdot\cdot+\left(  s+2\left(  m-1\right)  \pi\right)
^{n-1}\right]  ds\\
&  =\frac{1}{2\pi}\int_{0}^{2\pi}\varphi\left(  s\right)  \underbrace{\frac
{n}{m^{n}}\frac{1}{\left(  2\pi\right)  ^{n-1}}\left[  s^{n-1}+\left(
s+2\pi\right)  ^{n-1}+\left(  s+4\pi\right)  ^{n-1}+\cdot\cdot\cdot+\left(
s+2\left(  m-1\right)  \pi\right)  ^{n-1}\right]  }ds,
\end{align*}
where the underlined term$\ \underbrace{\cdot\cdot\cdot}$ is equal to
\begin{align*}
&  \underbrace{\cdot\cdot\cdot}\\
&  =\frac{n}{m^{n}}\left[  \left(  \frac{s}{2\pi}\right)  ^{n-1}+\left(
\frac{s}{2\pi}+1\right)  ^{n-1}+\left(  \frac{s}{2\pi}+2\right)  ^{n-1}%
+\cdot\cdot\cdot+\left(  \frac{s}{2\pi}+\left(  m-1\right)  \right)
^{n-1}\right]
\end{align*}
and for $s\in\left[  0,2\pi\right]  $ we have%
\begin{align}
&  \frac{n}{m^{n}}\left(  0^{n-1}+1^{n-1}+2^{n-1}+\cdot\cdot\cdot+\left(
m-1\right)  ^{n-1}\right) \nonumber\\
&  \leq\underbrace{\cdot\cdot\cdot}\leq\frac{n}{m^{n}}\left(  1^{n-1}%
+2^{n-1}+3^{n-1}+\cdot\cdot\cdot+m^{n-1}\right)  .
\end{align}
By the inequality
\begin{align}
&  \frac{n}{m^{n}}\left(  0^{n-1}+1^{n-1}+2^{n-1}+\cdot\cdot\cdot+\left(
m-1\right)  ^{n-1}\right) \nonumber\\
&  <\frac{n}{m^{n}}\int_{0}^{m}x^{n-1}dx=1<\frac{n}{m^{n}}\left(
1^{n-1}+2^{n-1}+3^{n-1}+\cdot\cdot\cdot+m^{n-1}\right)  \label{bet-2}%
\end{align}
with
\begin{align*}
&  \lim_{m\rightarrow\infty}\left[  \frac{n}{m^{n}}\left(  1^{n-1}%
+2^{n-1}+3^{n-1}+\cdot\cdot\cdot+m^{n-1}\right)  -\frac{n}{m^{n}}\left(
0^{n-1}+1^{n-1}+2^{n-1}+\cdot\cdot\cdot+\left(  m-1\right)  ^{n-1}\right)
\right] \\
&  =\lim_{m\rightarrow\infty}\left(  \frac{n}{m^{n}}m^{n-1}\right)  =0,
\end{align*}
we must have
\[
\left\{
\begin{array}
[c]{l}%
\lim_{m\rightarrow\infty}\frac{n}{m^{n}}\left(  0^{n-1}+1^{n-1}+2^{n-1}%
+\cdot\cdot\cdot+\left(  m-1\right)  ^{n-1}\right)  =1%
%TCIMACRO{\TeXButton{vspace}{\vspace{3mm}}}%
%BeginExpansion
\vspace{3mm}%
%EndExpansion
\\
\lim_{m\rightarrow\infty}\frac{n}{m^{n}}\left(  1^{n-1}+2^{n-1}+3^{n-1}%
+\cdot\cdot\cdot+m^{n-1}\right)  =1
\end{array}
\right.
\]
and hence $\lim_{m\rightarrow\infty}\underbrace{\cdot\cdot\cdot}=1,$ which
implies%
\begin{equation}
\lim_{m\rightarrow\infty}\frac{n}{\left(  2m\pi\right)  ^{n}}\int_{0}^{2m\pi
}\varphi\left(  r\right)  r^{n-1}dr=\frac{1}{2\pi}\int_{0}^{2\pi}%
\varphi\left(  s\right)  \left(  \lim_{m\rightarrow\infty}\underbrace
{\cdot\cdot\cdot}\right)  ds=\frac{1}{2\pi}\int_{0}^{2\pi}\varphi\left(
s\right)  ds. \label{conv-ave}%
\end{equation}
Finally we have%
\[
\lim_{m\rightarrow\infty}I=\lim_{m\rightarrow\infty}\left(  \frac{\left(
2m\pi\right)  ^{n}}{\left(  2m\pi+R\right)  ^{n}}\frac{n}{\left(
2m\pi\right)  ^{n}}\int_{0}^{2m\pi}\varphi\left(  r\right)  r^{n-1}dr\right)
=\frac{1}{2\pi}\int_{0}^{2\pi}\varphi\left(  s\right)  ds.
\]
The proof is done.$%
%TCIMACRO{\TeXButton{hfill}{\hfill}}%
%BeginExpansion
\hfill
%EndExpansion
\square$

\ \ \ \ 

As we have said in the paragraph before Theorem \ref{thm-main-2}, if
$\varphi\left(  \tau\right)  \ $has a slow oscillation\textbf{ }as
$\tau\rightarrow\infty\ $with $r<s,$ then it has the chance of passing to the
function $u\left(  0,t\right)  $. The following gives a simple way to find a
slow oscillation$\ $function$\ \varphi\left(  \tau\right)  $ on $[0,\infty
)\ $which will cause $\alpha<\beta.$

\begin{lemma}
\label{lem-slow}Let $g\left(  \tau\right)  $ be a $2\pi$-periodic non-constant
$C^{2}$ function defined on $[0,\infty)$ and let $G\left(  \tau\right)
=g\left(  \log\left(  \tau+1\right)  \right)  ,\ \tau\in\lbrack0,\infty).$
Then the radial function$\ \varphi\left(  \tau\right)  =\frac{\tau}%
{n}G^{\prime}\left(  \tau\right)  +G\left(  \tau\right)  ,\ \tau\in
\lbrack0,\infty),$ has a slow oscillation as $\tau\rightarrow\infty,$ lies in
the space $C^{0}[0,\infty)\bigcap L^{\infty}[0,\infty),$ and its average
integral $H\left(  \tau\right)  ,\ \tau\in\lbrack0,\infty),\ $%
satisfies$\ H\left(  \tau\right)  =g\left(  \log\left(  \tau+1\right)
\right)  ,\ \tau\in\lbrack0,\infty),\ $with%
\begin{equation}
p=\liminf_{\tau\rightarrow\infty}g\left(  \tau\right)  <\limsup_{\tau
\rightarrow\infty}g\left(  \tau\right)  =q, \label{lim}%
\end{equation}
which, by Theorem \ref{thm-K}, will also imply $\alpha<\beta.$
\end{lemma}

%

%TCIMACRO{\TeXButton{Proof}{\proof}}%
%BeginExpansion
\proof
%EndExpansion
By%
\[
\left\{
\begin{array}
[c]{l}%
\varphi\left(  \tau\right)  =\frac{\tau}{n}\frac{g^{\prime}\left(  \log\left(
\tau+1\right)  \right)  }{\tau+1}+g\left(  \log\left(  \tau+1\right)  \right)
,\ \ \ \varphi\left(  0\right)  =g\left(  0\right)  ,%
%TCIMACRO{\TeXButton{vspace}{\vspace{3mm}}}%
%BeginExpansion
\vspace{3mm}%
%EndExpansion
\\
\varphi^{\prime}\left(  \tau\right)  =\left(  1+\frac{1}{n}\right)
\frac{g^{\prime}\left(  \log\left(  \tau+1\right)  \right)  }{\tau+1}%
+\frac{\tau\left[  g^{\prime\prime}\left(  \log\left(  \tau+1\right)  \right)
-g^{\prime}\left(  \log\left(  \tau+1\right)  \right)  \right]  }{n\left(
\tau+1\right)  ^{2}},\ \ \ \tau\in\lbrack0,\infty),
\end{array}
\right.
\]
we see that$\ \varphi\ $is a slow oscillation function on $[0,\infty),\ $lying
in the space $C^{0}[0,\infty)\bigcap L^{\infty}[0,\infty).\ $Now%
\begin{align*}
H\left(  \tau\right)   &  =\frac{n}{\tau^{n}}\int_{0}^{\tau}\varphi\left(
r\right)  r^{n-1}dr=\frac{n}{\tau^{n}}\int_{0}^{\tau}\left(  \frac{r}%
{n}G^{\prime}\left(  r\right)  +G\left(  r\right)  \right)  r^{n-1}dr\\
&  =\frac{1}{\tau^{n}}\int_{0}^{\tau}\frac{d}{dr}\left(  r^{n}G\left(
r\right)  \right)  dr=G\left(  \tau\right)  =g\left(  \log\left(
\tau+1\right)  \right)  ,\ \ \ \tau\in\left(  0,\infty\right)
\end{align*}
and$\ H\left(  0\right)  =G\left(  0\right)  =g\left(  0\right)
.\ $Therefore, the inequality (\ref{lim}) follows.$%
%TCIMACRO{\TeXButton{hfill}{\hfill}}%
%BeginExpansion
\hfill
%EndExpansion
\square$

\ \ \ 

With the help of Lemma \ref{lem-ave}, we are ready to prove Theorem
\ref{thm-main-2}.\ For any $r<\alpha<\beta<s,$ we divide the proof into two
cases.$%
%TCIMACRO{\TeXButton{vspace}{\vspace{3mm}}}%
%BeginExpansion
\vspace{3mm}%
%EndExpansion
$

\textbf{Case 1:\ }$r+s=\alpha+\beta,\ r<\alpha<\beta<s.%
%TCIMACRO{\TeXButton{vspace}{\vspace{3mm}}}%
%BeginExpansion
\vspace{3mm}%
%EndExpansion
$

In this case we can use the representation formula (\ref{CONV-rad}) and,
similar to (\ref{H-tau-m}), choose
\begin{equation}
\varphi\left(  \tau\right)  =\frac{s-r}{2}\sin\left(  m\log\left(
\tau+1\right)  \right)  +\frac{s+r}{2},\ \ \ \varphi\left(  0\right)
=\frac{s+r}{2},\ \ \ \tau\in\lbrack0,\infty), \label{phi-4}%
\end{equation}
and perform the same argument as in the proof of Theorem \ref{thm-main-1} to
prescribe them.\ More precisely, here we need to choose $m\in\left(
0,\infty\right)  $ to satisfy $\sqrt{\tilde{A}^{2}\left(  m\right)  +\tilde
{B}^{2}\left(  m\right)  }=\frac{\beta-\alpha}{s-r}\in\left(  0,1\right)  ,$
where now%
\begin{equation}
\left\{
\begin{array}
[c]{l}%
\tilde{A}\left(  m\right)  =\dfrac{n\omega\left(  n\right)  }{\pi^{n/2}}%
\int_{0}^{\infty}e^{-z^{2}}z^{n-1}\cos\left(  m\log z\right)  dz\in\left(
-1,1\right)  ,%
%TCIMACRO{\TeXButton{vspace}{\vspace{3mm}}}%
%BeginExpansion
\vspace{3mm}%
%EndExpansion
\\
\tilde{B}\left(  m\right)  =\dfrac{n\omega\left(  n\right)  }{\pi^{n/2}}%
\int_{0}^{\infty}e^{-z^{2}}z^{n-1}\sin\left(  m\log z\right)  dz\in\left(
-1,1\right)  ,
\end{array}
\right.  \label{AmBm-2}%
\end{equation}
with
\begin{equation}
0<\tilde{A}^{2}\left(  m\right)  +\tilde{B}^{2}\left(  m\right)
<\frac{n\omega\left(  n\right)  }{\pi^{n/2}}\int_{0}^{\infty}e^{-z^{2}}%
z^{n-1}\left[  \cos^{2}\left(  m\log z\right)  +\sin^{2}\left(  m\log
z\right)  \right]  dz=1.
\end{equation}
Note that in this case the chosen function $\varphi\left(  \tau\right)  ,$
like $H\left(  \tau\right)  \ $in (\ref{H-tau-m}), has a \textbf{slow
oscillation} on $[0,\infty).$

$%
%TCIMACRO{\TeXButton{vspace}{\vspace{3mm}}}%
%BeginExpansion
\vspace{3mm}%
%EndExpansion
$

\textbf{Case 2:}$\ r+s\neq\alpha+\beta,\ r<\alpha<\beta<s.%
%TCIMACRO{\TeXButton{vspace}{\vspace{3mm}}}%
%BeginExpansion
\vspace{3mm}%
%EndExpansion
$

In this case, we may assume $r+s>\alpha+\beta\ $(the treatment for the case
$r+s<\alpha+\beta$ is similar). We will choose$\ \varphi\left(  \tau\right)
\in C^{0}[0,\infty)\bigcap L^{\infty}[0,\infty)$ to be the \textbf{sum} of two
functions $\varphi_{1}\left(  \tau\right)  +\varphi_{2}\left(  \tau\right)  ,$
where $\varphi_{1}\left(  \tau\right)  $ has a \textbf{slow oscillation}, but
$\varphi_{2}\left(  \tau\right)  $ has a \textbf{regular oscillation} (which
is $2\pi$\textbf{-periodic}).$\ $

Let $\lambda\ $be the number satisfying $r+\lambda=\alpha+\beta.\ $Since
$r<\alpha<\beta<s$ and $r+s>\alpha+\beta,$ we have$\ \beta<\lambda<s.$ We can
write $\lambda$ as $\lambda=\varepsilon+\delta$ for some small $\varepsilon>0$
and $\delta\in\mathbb{R}\ $so that we have%
\begin{equation}
r+\varepsilon<\alpha<\beta<\delta<s.
\end{equation}
Now we have $\left(  r+\varepsilon\right)  +\delta=\alpha+\beta\ $and by Case
1, we can find a radial $\varphi_{1}\in C^{0}[0,\infty)\bigcap L^{\infty
}[0,\infty)\ $such that it and its convolution solution $u_{1}\left(
x,t\right)  $ satisfy%
\begin{equation}
\liminf_{\tau\rightarrow\infty}\varphi_{1}\left(  \tau\right)  =r+\varepsilon
<\liminf_{t\rightarrow\infty}u_{1}\left(  0,t\right)  =\alpha<\limsup
_{t\rightarrow\infty}u_{1}\left(  0,t\right)  =\beta<\limsup_{\tau
\rightarrow\infty}\varphi_{1}\left(  \tau\right)  =\delta. \label{cond1}%
\end{equation}
Next let $\varphi_{2}\left(  \tau\right)  \ $be a $2\pi$-periodic continuous
function on$\ \tau\in\lbrack0,\infty)$ which satisfies
\begin{equation}
\frac{1}{2\pi}\int_{0}^{2\pi}\varphi_{2}\left(  \tau\right)  d\tau
=0,\ \ \ \min_{\tau\in\left[  0,2\pi\right]  }\varphi_{2}\left(  \tau\right)
=-\varepsilon<0,\ \ \ \max_{\tau\in\left[  0,2\pi\right]  }\varphi_{2}\left(
\tau\right)  =s-\delta>0. \label{cond2}%
\end{equation}
Such a function $\varphi_{2}\ $clearly exists and, together with its
convolution solution $u_{2}\left(  x,t\right)  ,\ $they satisfy%
\begin{equation}
\liminf_{\tau\rightarrow\infty}\varphi_{2}\left(  \tau\right)  =-\varepsilon
<\lim_{t\rightarrow\infty}u_{2}\left(  0,t\right)  =\frac{1}{2\pi}\int
_{0}^{2\pi}\varphi_{2}\left(  \tau\right)  d\tau=0<\limsup_{\tau
\rightarrow\infty}\varphi_{2}\left(  \tau\right)  =s-\delta, \label{cond3}%
\end{equation}
due to the result in Lemma \ref{lem-ave}. Finally, we set $\varphi\left(
\tau\right)  =\varphi_{1}\left(  \tau\right)  +\varphi_{2}\left(  \tau\right)
,\ \tau\in\lbrack0,\infty).\ $We have$\ \varphi\in C^{0}[0,\infty)\bigcap
L^{\infty}[0,\infty)\ $and its convolution solution $u\left(  x,t\right)  $
satisfy$\ u\left(  x,t\right)  =u_{1}\left(  x,t\right)  +u_{2}\left(
x,t\right)  \ $for all $\left(  x,t\right)  \in\mathbb{R}^{n}\times
\lbrack0,\infty).\ $Since$\ u_{2}\left(  0,t\right)  \rightarrow0\ $as
$t\rightarrow\infty,$ we have%
\begin{equation}
\left\{
\begin{array}
[c]{l}%
\liminf_{t\rightarrow\infty}u\left(  0,t\right)  =\liminf_{t\rightarrow\infty
}u_{1}\left(  0,t\right)  =\alpha%
%TCIMACRO{\TeXButton{vspace}{\vspace{3mm}}}%
%BeginExpansion
\vspace{3mm}%
%EndExpansion
\\
\limsup_{t\rightarrow\infty}u\left(  0,t\right)  =\limsup_{t\rightarrow\infty
}u_{1}\left(  0,t\right)  =\beta.
\end{array}
\right.  \label{ab}%
\end{equation}
It remains to look at $\varphi\left(  \tau\right)  $ for $\tau\rightarrow
\infty.\ $By (\ref{phi-4}) we know $\varphi_{1}\left(  \tau\right)  $ is a
\textbf{slow oscillation} function on $[0,\infty)$ with the form%
\begin{equation}
\varphi_{1}\left(  \tau\right)  =\frac{\delta-\left(  r+\varepsilon\right)
}{2}\sin\left(  m\log\left(  \tau+1\right)  \right)  +\frac{\delta+\left(
r+\varepsilon\right)  }{2},\ \ \ \varphi_{1}\left(  0\right)  =\frac
{\delta+\left(  r+\varepsilon\right)  }{2},\ \ \ \tau\in\lbrack0,\infty)
\label{phi-5}%
\end{equation}
for some fixed$\ m\in\left(  0,\infty\right)  \ $satisfying$\ \sqrt{\tilde
{A}^{2}\left(  m\right)  +\tilde{B}^{2}\left(  m\right)  }=\frac{\beta-\alpha
}{\delta-\left(  r+\varepsilon\right)  }\in\left(  0,1\right)  .\ $On the
other hand, $\varphi_{2}\left(  \tau\right)  $ is a \textbf{regular
oscillation} function ($2\pi$\textbf{-periodic function})\ on $[0,\infty),$
which oscillates between $-\varepsilon\ $and $s-\delta.\ $It is not difficult
to see that there exists a sequence $\tau_{j}\rightarrow\infty\ $so that%
\begin{equation}
\sin\left(  m\log\left(  \tau_{j}+1\right)  \right)  \rightarrow
-1\ \ \ \text{and\ \ \ }\varphi_{2}\left(  \tau_{j}\right)  \rightarrow
-\varepsilon=\min_{\tau\in\left[  0,2\pi\right]  }\varphi_{2}\left(
\tau\right)  <0 \label{sin-1}%
\end{equation}
as $j\rightarrow\infty.\ $ Similarly, there exists a sequence $\tau
_{k}\rightarrow\infty\ $so that%
\begin{equation}
\sin\left(  m\log\left(  \tau_{k}+1\right)  \right)  \rightarrow
1\ \ \ \text{and\ \ \ }\varphi_{2}\left(  \tau_{k}\right)  \rightarrow
s-\delta=\max_{\tau\in\left[  0,2\pi\right]  }\varphi_{2}\left(  \tau\right)
>0 \label{sin1}%
\end{equation}
as $k\rightarrow\infty.$ By (\ref{phi-5}), (\ref{sin1}),\ and (\ref{sin-1}),
we have%
\begin{align}
\liminf_{\tau\rightarrow\infty}\varphi\left(  \tau\right)   &  =\liminf
_{\tau\rightarrow\infty}\left(  \varphi_{1}\left(  \tau\right)  +\varphi
_{2}\left(  \tau\right)  \right) \nonumber\\
&  =\liminf_{\tau\rightarrow\infty}\varphi_{1}\left(  \tau\right)
+\liminf_{\tau\rightarrow\infty}\varphi_{2}\left(  \tau\right)  =\left(
r+\varepsilon\right)  -\varepsilon=r \label{rr}%
\end{align}
and%
\begin{align}
\limsup_{\tau\rightarrow\infty}\varphi\left(  \tau\right)   &  =\limsup
_{\tau\rightarrow\infty}\left(  \varphi_{1}\left(  \tau\right)  +\varphi
_{2}\left(  \tau\right)  \right) \nonumber\\
&  =\limsup_{\tau\rightarrow\infty}\varphi_{1}\left(  \tau\right)
+\limsup_{\tau\rightarrow\infty}\varphi_{2}\left(  \tau\right)  =\delta
+\left(  s-\delta\right)  =s. \label{ss}%
\end{align}
The proof is done due to (\ref{ab}), (\ref{rr}), and\ (\ref{ss}).

Combining Case 1 and Case 2,\ the proof of Theorem \ref{thm-main-2} is now
complete.$%
%TCIMACRO{\TeXButton{hfill}{\hfill}}%
%BeginExpansion
\hfill
%EndExpansion
\square${}

\begin{remark}
This is to explain the existence of a sequence $\tau_{k}\rightarrow\infty\ $so
that (\ref{sin1}) holds.\ We first note that, for given $\varepsilon
>0,\ $there\ exists$\ \delta>0\ $(depending only on $\varepsilon$), such that
if $m\log\left(  \tau+1\right)  $ (here $m>0$ is fixed)\ lies in the interval
$\left(  2k\pi+\frac{\pi}{2}-\delta,2k\pi+\frac{\pi}{2}+\delta\right)  \ $for
\textbf{some} $k\in\mathbb{N},$ then we have$\ \left\vert \sin\left(
m\log\left(  \tau+1\right)  \right)  -1\right\vert <\varepsilon,$ which is
equivalent\ for $\tau+1$ to lie in the interval
\begin{equation}
\tau+1\in\left(  \exp\left(  \frac{2k\pi+\frac{\pi}{2}-\delta}{m}\right)
,\ \exp\left(  \frac{2k\pi+\frac{\pi}{2}+\delta}{m}\right)  \right)  .
\label{sin}%
\end{equation}
We can\ choose $k\in\mathbb{N}$ sufficiently large so that the length of the
interval in (\ref{sin}) is sufficiently large.$\ $By this observation, we
clearly have (\ref{sin1}). The reason for (\ref{sin-1}) to hold is similar.
\end{remark}

\subsection{The remaining cases not covered by Theorem \ref{thm-main-2}.}

For prescribing general$\ r\leq\alpha\leq\beta\leq s,\ $the remaining cases
not covered by Theorem \ref{thm-main-2} are discussed in the following:

\begin{theorem}
\label{thm-main-3}For any four arbitrary finite numbers $r,\ \alpha
,\ \beta,\ s\ $in one of the following cases:
\begin{equation}
\left\{
\begin{array}
[c]{l}%
\left(  1\right)  .\ r=\alpha<\beta<s,\ \ \ \ \ \left(  2\right)
.\ r<\alpha=\beta<s,\ \ \ \ \ \left(  3\right)  .\ r<\alpha<\beta=s,%
%TCIMACRO{\TeXButton{vspace}{\vspace{3mm}}}%
%BeginExpansion
\vspace{3mm}%
%EndExpansion
\\
\left(  4\right)  .\ r=\alpha=\beta<s,\ \ \ \ \ \left(  5\right)
.\ r=\alpha<\beta=s,\ \ \ \ \ \left(  6\right)  .\ r<\alpha=\beta=s,
\end{array}
\right.  \label{6}%
\end{equation}
one can choose a\ suitable radial initial data $\varphi\left(  x\right)
=\varphi\left(  \left\vert x\right\vert \right)  \in C^{0}[0,\infty)\bigcap
L^{\infty}[0,\infty)\ $satisfying
\begin{equation}
\liminf_{\tau\rightarrow\infty}\varphi\left(  \tau\right)  =r,\ \ \ \limsup
_{\tau\rightarrow\infty}\varphi\left(  \tau\right)  =s
\end{equation}
and its corresponding convolution solution $u\left(  x,t\right)  $ in
(\ref{CONV})\ satisfies%
\begin{equation}
\liminf_{t\rightarrow\infty}u\left(  0,t\right)  =\alpha,\ \ \ \limsup
_{t\rightarrow\infty}u\left(  0,t\right)  =\beta.
\end{equation}

\end{theorem}

\begin{remark}
The case $\ r=\alpha=\beta=s\ $is trivial.\ Just choose $\varphi\left(
x\right)  $ to be a constant function.
\end{remark}

%

%TCIMACRO{\TeXButton{Proof}{\proof}}%
%BeginExpansion
\proof
%EndExpansion
$\left(  2\right)  .\ $For $r<\alpha=\beta<s,\ $we can choose $\varphi\left(
\tau\right)  ,\ \tau\in\lbrack0,\infty),\ $to be any $2\pi$-periodic function
with average value $\alpha=\beta,$ maximum value $s,\ $minimum value $r.\ $By
Lemma \ref{lem-ave}, it can be achieved.\ $%
%TCIMACRO{\TeXButton{vspace}{\vspace{3mm}}}%
%BeginExpansion
\vspace{3mm}%
%EndExpansion
$

$\left(  5\right)  .\ $For $r=\alpha<\beta=s,\ $we choose $\varphi\left(
\tau\right)  $ to be the\textbf{ extremely slow oscillation function}
\begin{equation}
\varphi\left(  \tau\right)  =\frac{\beta-\alpha}{2}\sin\left[  \log\left(
\log\left(  \tau+2\right)  \right)  \right]  +\frac{\beta+\alpha}%
{2},\ \ \ \tau\in\lbrack0,\infty), \label{ini-data}%
\end{equation}
which has%
\begin{equation}
\liminf_{\tau\rightarrow\infty}\varphi\left(  \tau\right)  =\alpha
,\ \ \ \limsup_{\tau\rightarrow\infty}\varphi\left(  \tau\right)  =\beta.
\label{phi-conv}%
\end{equation}
By (\ref{CONV-rad}), its corresponding $u\left(  0,t\right)  $ is given by%
\begin{equation}
u\left(  0,t\right)  =\frac{\beta-\alpha}{2}\dfrac{n\omega\left(  n\right)
}{\pi^{n/2}}\int_{0}^{\infty}e^{-z^{2}}z^{n-1}\sin\left[  \log\left(
\log\left(  \sqrt{4t}z+2\right)  \right)  \right]  dz+\frac{\beta+\alpha}{2},
\label{u0t}%
\end{equation}
where one can write%
\[
\sin\left[  \log\left(  \log\left(  \sqrt{4t}z+2\right)  \right)  \right]
=\sin\left[  \log\left(  \log\sqrt{4t}\right)  +g\left(  t,z\right)  \right]
,
\]
with
\[
g\left(  t,z\right)  =\log\left(  1+\frac{1}{\log\sqrt{4t}}\log\left(
z+\frac{2}{\sqrt{4t}}\right)  \right)  ,\ \ \ t,\ z\in\left(  0,\infty\right)
.
\]
We see that $\lim_{t\rightarrow\infty}g\left(  t,z\right)  =0\ $for fixed
$z\in\left(  0,\infty\right)  .$ Hence the Lebesgue Dominated Convergence
Theorem implies%
\begin{equation}
\lim_{t\rightarrow\infty}\left\vert u\left(  0,t\right)  -\left[  \frac
{\beta-\alpha}{2}\sin\left(  \log\left(  \log\sqrt{4t}\right)  \right)
+\frac{\beta+\alpha}{2}\right]  \right\vert =0 \label{u0t-conv}%
\end{equation}
and so
\begin{equation}
\liminf_{t\rightarrow\infty}u\left(  0,t\right)  =\alpha,\ \ \ \limsup
_{t\rightarrow\infty}u\left(  0,t\right)  =\beta. \label{u-conv}%
\end{equation}

$\left(  4\right)  .\ $For$\ r=\alpha=\beta<s,$ we\ choose $\varphi\in
C^{0}[0,\infty)\bigcap L^{\infty}[0,\infty)$ with $\liminf_{\tau
\rightarrow\infty}\varphi\left(  \tau\right)  =r=\alpha,\ \limsup
_{\tau\rightarrow\infty}\varphi\left(  \tau\right)  =s,$ and for most $\tau
\in\lbrack0,\infty)\ $the value of $\varphi\left(  \tau\right)  $ is equal to
$\alpha\ $and for $\varphi\left(  \tau\right)  $ not equal to $\alpha,$ it
looks like a thin bump with height $s-\alpha\ $and the supports of these bumps
are spreading further and further apart. By this choice of $\varphi$ we have%
\[
\lim_{\tau\rightarrow\infty}H\left(  \tau\right)  =\lim_{\tau\rightarrow
\infty}\frac{n}{\tau^{n}}\int_{0}^{\tau}\varphi\left(  r\right)
r^{n-1}dr=\alpha,
\]
which implies $\lim_{t\rightarrow\infty}u\left(  0,t\right)  =\alpha\ $due to
Theorem \ref{thm-K}.\ Hence, the case $r=\alpha=\beta<s$ is achieved.\ $%
%TCIMACRO{\TeXButton{vspace}{\vspace{3mm}}}%
%BeginExpansion
\vspace{3mm}%
%EndExpansion
$

$\left(  6\right)  .\ $For$\ r<\alpha=\beta=s,\ $the construction of $\varphi$
is similar to that in $\left(  4\right)  \ $except that we reverse the role of
$r\ $and $s.$ $%
%TCIMACRO{\TeXButton{vspace}{\vspace{3mm}}}%
%BeginExpansion
\vspace{3mm}%
%EndExpansion
$

$\left(  1\right)  .\ $For$\ r=\alpha<\beta<s,\ $we first choose $\varphi
_{1}\left(  \tau\right)  $ to be the function given by (\ref{ini-data}),
which, together with its corresponding $u_{1}\left(  0,t\right)  ,$ will
satisfy (\ref{phi-conv}) and (\ref{u-conv}). Next, we choose $\varphi
_{2}\left(  \tau\right)  .\ $Let $\varepsilon>0\ $be a fixed small number and
for each$\ m\in\mathbb{N},\ $let
\[
\tau_{m}=\exp\left(  \exp\left(  2m\pi+\frac{\pi}{2}\right)  \right)
-2,\ \ \ \ \tilde{\tau}_{m}=\exp\left(  \exp\left(  2m\pi+\frac{3\pi}%
{2}\right)  \right)  -2,
\]
where we note that $\varphi_{1}\left(  \tilde{\tau}_{m}\right)  =\alpha
,\ \varphi_{1}\left(  \tau_{m}\right)  =\beta.\ $We require$\ \varphi
_{2}\left(  \tau\right)  $ to be a nonnegative function satisfying%
\[
\varphi_{2}\left(  \tau\right)  =\left\{
\begin{array}
[c]{l}%
\frac{s-\beta}{\varepsilon}\left(  \tau-\left(  \tau_{m}-\varepsilon\right)
\right)  ,\ \ \ \tau\in\left[  \tau_{m}-\varepsilon,\tau_{m}\right]
%TCIMACRO{\TeXButton{vspace}{\vspace{3mm}}}%
%BeginExpansion
\vspace{3mm}%
%EndExpansion
\\
-\frac{s-\beta}{\varepsilon}\left(  \tau-\left(  \tau_{m}+\varepsilon\right)
\right)  ,\ \ \ \tau\in\left[  \tau_{m},\tau_{m}+\varepsilon\right]
%TCIMACRO{\TeXButton{vspace}{\vspace{3mm}}}%
%BeginExpansion
\vspace{3mm}%
%EndExpansion
\\
0,\ \ \ \text{otherwise,\ \ \ }\tau\in\lbrack0,\infty).
\end{array}
\right.
\]
It satisfies $\liminf_{\tau\rightarrow\infty}\varphi_{2}\left(  \tau\right)
=0\ $and$\ \limsup_{\tau\rightarrow\infty}\varphi_{2}\left(  \tau\right)
=s-\beta\ $and%
\[
\lim_{\tau\rightarrow\infty}H_{2}\left(  \tau\right)  =\lim_{\tau
\rightarrow\infty}\frac{n}{\tau^{n}}\int_{0}^{\tau}\varphi_{2}\left(
r\right)  r^{n-1}dr=0.
\]
Hence its corresponding $u_{2}\left(  0,t\right)  $ satisfies$\ \lim
_{t\rightarrow\infty}u_{2}\left(  0,t\right)  =0\ $due to Theorem \ref{thm-K}.
Now set $\varphi=\varphi_{1}+\varphi_{2}\in C^{0}[0,\infty)\bigcap L^{\infty
}[0,\infty).\ $We will have $u\left(  0,t\right)  =u_{1}\left(  0,t\right)
+u_{2}\left(  0,t\right)  $ with$\ $
\[
\left\{
\begin{array}
[c]{l}%
\liminf_{t\rightarrow\infty}u\left(  0,t\right)  =\liminf_{t\rightarrow\infty
}u_{1}\left(  0,t\right)  =\alpha%
%TCIMACRO{\TeXButton{vspace}{\vspace{3mm}}}%
%BeginExpansion
\vspace{3mm}%
%EndExpansion
\\
\limsup_{t\rightarrow\infty}u\left(  0,t\right)  =\limsup_{t\rightarrow\infty
}u_{1}\left(  0,t\right)  =\beta.
\end{array}
\right.
\]
Also, by$\ \varphi\left(  \tilde{\tau}_{m}\right)  =\varphi_{1}\left(
\tilde{\tau}_{m}\right)  +\varphi_{2}\left(  \tilde{\tau}_{m}\right)
=\alpha+0=\alpha\ $and$\ \varphi\left(  \tau_{m}\right)  =\varphi_{1}\left(
\tau_{m}\right)  +\varphi_{2}\left(  \tau_{m}\right)  =\beta+\left(
s-\beta\right)  =s,$ we have
\[
\left\{
\begin{array}
[c]{l}%
\liminf_{\tau\rightarrow\infty}\varphi\left(  \tau\right)  =\liminf
_{\tau\rightarrow\infty}\left(  \varphi_{1}\left(  \tau\right)  +\varphi
_{2}\left(  \tau\right)  \right)  =\alpha%
%TCIMACRO{\TeXButton{vspace}{\vspace{3mm}}}%
%BeginExpansion
\vspace{3mm}%
%EndExpansion
\\
\limsup_{\tau\rightarrow\infty}\varphi\left(  \tau\right)  =\limsup
_{\tau\rightarrow\infty}\left(  \varphi_{1}\left(  \tau\right)  +\varphi
_{2}\left(  \tau\right)  \right)  =\beta+\left(  s-\beta\right)  =s.
\end{array}
\right.
\]
$\ $ $\ $ \ $\ $

$\left(  3\right)  \ $For$\ r<\alpha<\beta=s,\ $the construction of $\varphi$
is similar to that in $\left(  1\right)  .%
%TCIMACRO{\TeXButton{vspace}{\vspace{3mm}}}%
%BeginExpansion
\vspace{3mm}%
%EndExpansion
$

The proof of Theorem \ref{thm-main-3} is now complete.$%
%TCIMACRO{\TeXButton{hfill}{\hfill}}%
%BeginExpansion
\hfill
%EndExpansion
\square$

\section{Some side issues.}

\subsection{An example of prescribing four different\ numbers not satisfying
(\ref{sym-cond}) in Theorem \ref{thm-main-1}.\label{sec-exam}}

This is related to Theorem \ref{thm-main-1}.\ Until now, we still do not know
how to \textbf{prescribe} the oscillation of $H\left(  \tau\right)  \ $and
$u\left(  0,t\right)  $ for four\textbf{ }arbitrary numbers $p<\alpha
<\beta<q\ $\textbf{not} satisfying the condition (\ref{sym-cond}). This will
be an interesting problem to explore. However, it is not difficult to
construct a particular example not satisfying (\ref{sym-cond}). We have:

\begin{lemma}
\label{lem-not}There exists a radial function $\varphi\in C^{0}[0,\infty
)\bigcap L^{\infty}[0,\infty)$ so that its average integral $H\left(
\tau\right)  \ $and the convolution solution (\ref{CONV}) of the problem
(\ref{heat})\ satisfy $p<\alpha<\beta<q$ with $p+q\neq\alpha+\beta.$
\end{lemma}

%

%TCIMACRO{\TeXButton{Proof}{\proof}}%
%BeginExpansion
\proof
%EndExpansion
The idea of breaking the symmetry, unlike (\ref{Hm-tau}), is to choose
$H\left(  \tau\right)  $ to be equal to the sum of \textbf{several }(at least
two)\textbf{\ }different\ slow oscillation functions.\ For simplicity, we only
look at the case $n=1\ $(the construction for $n>1\ $is similar, but we cannot
find numerical values involving a general variable $n$)\ and require $H\left(
\tau\right)  $ to be equal to the sum of\textbf{ two }slow oscillation
functions, given by
\begin{align}
H\left(  \tau\right)   &  =\frac{n}{\tau^{n}}\int_{0}^{\tau}\varphi\left(
r\right)  r^{n-1}dr=\frac{1}{\tau}\int_{0}^{\tau}\varphi\left(  r\right)
dr\ \ \ \left(  n=1\right) \nonumber\\
&  =\sin\left(  \log\left(  \tau+1\right)  \right)  +\sin\left(  2\log\left(
\tau+1\right)  \right)  ,\ \ \ \tau\in\left(  0,\infty\right)  ,\ \ \ H\left(
0\right)  =0.
\end{align}
Similar to (\ref{3-tau}), one can find the corresponding initial radial
function $\varphi\in C^{0}[0,\infty)\bigcap L^{\infty}[0,\infty).$ With the
help of Maple software, we can evaluate
\[
\left\{
\begin{array}
[c]{l}%
p=\liminf_{\tau\rightarrow\infty}H\left(  \tau\right)  =\min_{x\in\left[
0,2\pi\right]  }\left(  \sin x+\sin2x\right)  \approx-1.760172593%
%TCIMACRO{\TeXButton{vspace}{\vspace{3mm}}}%
%BeginExpansion
\vspace{3mm}%
%EndExpansion
\\
q=\limsup_{\tau\rightarrow\infty}H\left(  \tau\right)  =\max_{x\in\left[
0,2\pi\right]  }\left(  \sin x+\sin2x\right)  \approx1.760172593.
\end{array}
\right.
\]
By the representation formula $u\left(  0,t\right)  =\frac{4}{\sqrt{\pi}}%
\int_{0}^{\infty}e^{-z^{2}}z^{2}H\left(  \sqrt{4t}z\right)  dz,\ t\in\left(
0,\infty\right)  ,\ $we have
\[
u\left(  0,t\right)  =\left\{
\begin{array}
[c]{l}%
A\left(  t\right)  \sin\left(  \log\sqrt{4t}\right)  +B\left(  t\right)
\cos\left(  \log\sqrt{4t}\right)
%TCIMACRO{\TeXButton{vspace}{\vspace{3mm}}}%
%BeginExpansion
\vspace{3mm}%
%EndExpansion
\\
+C\left(  t\right)  \sin\left(  2\log\sqrt{4t}\right)  +D\left(  t\right)
\cos\left(  2\log\sqrt{4t}\right)  ,\ \ \ t\in\left(  0,\infty\right)  ,
\end{array}
\right.
\]
where
\[
\left\{
\begin{array}
[c]{l}%
A\left(  t\right)  =\frac{4}{\sqrt{\pi}}\int_{0}^{\infty}e^{-z^{2}}z^{2}%
\cos\left(  \log\left(  z+\frac{1}{\sqrt{4t}}\right)  \right)  dz%
%TCIMACRO{\TeXButton{vspace}{\vspace{3mm}}}%
%BeginExpansion
\vspace{3mm}%
%EndExpansion
\\
B\left(  t\right)  =\frac{4}{\sqrt{\pi}}\int_{0}^{\infty}e^{-z^{2}}z^{2}%
\sin\left(  \log\left(  z+\frac{1}{\sqrt{4t}}\right)  \right)  dz%
%TCIMACRO{\TeXButton{vspace}{\vspace{3mm}}}%
%BeginExpansion
\vspace{3mm}%
%EndExpansion
\\
C\left(  t\right)  =\frac{4}{\sqrt{\pi}}\int_{0}^{\infty}e^{-z^{2}}z^{2}%
\cos\left(  2\log\left(  z+\frac{1}{\sqrt{4t}}\right)  \right)  dz%
%TCIMACRO{\TeXButton{vspace}{\vspace{3mm}}}%
%BeginExpansion
\vspace{3mm}%
%EndExpansion
\\
D\left(  t\right)  =\frac{4}{\sqrt{\pi}}\int_{0}^{\infty}e^{-z^{2}}z^{2}%
\sin\left(  2\log\left(  z+\frac{1}{\sqrt{4t}}\right)  \right)  dz,\ \ \ t\in
\left(  0,\infty\right)  .
\end{array}
\right.
\]
As $t\rightarrow\infty,$ we have
\[
\left\{
\begin{array}
[c]{l}%
\lim_{t\rightarrow\infty}A\left(  t\right)  =A\approx0.892253317,\ \ \ \lim
_{t\rightarrow\infty}B\left(  t\right)  =B\approx0.030945895%
%TCIMACRO{\TeXButton{vspace}{\vspace{3mm}}}%
%BeginExpansion
\vspace{3mm}%
%EndExpansion
\\
\lim_{t\rightarrow\infty}C\left(  t\right)  =C\approx0.649173672,\ \ \ \lim
_{t\rightarrow\infty}D\left(  t\right)  =D\approx0.099535090,
\end{array}
\right.
\]
and so
\begin{align*}
&  \lim_{t\rightarrow\infty}\left\vert u\left(  0,t\right)  -\left[
A\sin\left(  \log\sqrt{4t}\right)  +B\cos\left(  \log\sqrt{4t}\right)
+C\sin\left(  2\log\sqrt{4t}\right)  +D\cos\left(  2\log\sqrt{4t}\right)
\right]  \right\vert \\
&  =0,
\end{align*}
which gives%
\[
\left\{
\begin{array}
[c]{l}%
\alpha=\liminf_{t\rightarrow\infty}u\left(  0,t\right)  =\min_{x\in\left[
0,2\pi\right]  }\left(  A\sin x+B\cos x+C\sin2x+D\cos2x\right)  \approx
-1.369211837%
%TCIMACRO{\TeXButton{vspace}{\vspace{3mm}}}%
%BeginExpansion
\vspace{3mm}%
%EndExpansion
\\
\beta=\limsup_{t\rightarrow\infty}u\left(  0,t\right)  =\max_{x\in\left[
0,2\pi\right]  }\left(  A\sin x+B\cos x+C\sin2x+D\cos2x\right)  \approx
1.328017886
\end{array}
\right.
\]
and we find that $p+q\neq\alpha+\beta.\
%TCIMACRO{\TeXButton{hfill}{\hfill}}%
%BeginExpansion
\hfill
%EndExpansion
\square$

\begin{remark}
We have used different\ software to evaluate the above numerical values of
$p,\ \alpha,\ \beta,\ q\ $and obtain\ the same values.\ 
\end{remark}

\subsection{The behavior of $u\left(  x,t\right)  $ as $\left\vert
x\right\vert \rightarrow\infty.$}

In Theorem \ref{thm-main-1} and Theorem \ref{thm-main-2}\ we only look at the
behavior of $u\left(  x,t\right)  $ for fixed $x\ $with $t\rightarrow
\infty.\ $We may as well look at the behavior of $u\left(  x,t\right)  $ for
fixed $t\ $with $\left\vert x\right\vert \rightarrow\infty.\ $In the following
lemma, we discuss such property for a particular radial initial data.\ We
choose $\varphi\left(  x\right)  =\varphi\left(  \left\vert x\right\vert
\right)  \in C^{0}[0,\infty)\bigcap L^{\infty}[0,\infty)$ to be the function
in (\ref{3-tau}) with $m=1.\ $Its average integral $H\left(  \tau\right)
\ $is equal to $\sin\left(  \log\left(  \tau+1\right)  \right)  ,$ $\tau
\in\lbrack0,\infty),\ $which is our typical slow oscillation function. $\ $

\begin{lemma}
Let $u\left(  x,t\right)  $ be the solution of the problem\ (\ref{heat}) given
by (\ref{CONV}) with $u\left(  x,0\right)  =\varphi\left(  \left\vert
x\right\vert \right)  ,$ $x\in\mathbb{R}^{n},\ $where$\ \varphi\left(
\tau\right)  $ is given by%
\begin{equation}
\varphi\left(  \tau\right)  =\sin\left(  \log\left(  \tau+1\right)  \right)
+\frac{\tau}{n\left(  \tau+1\right)  }\cos\left(  \log\left(  \tau+1\right)
\right)  ,\ \ \ \tau\in\lbrack0,\infty). \label{phi-3}%
\end{equation}
Then for each fixed $t>0$ we have%
\begin{equation}
\lim_{\left\vert x\right\vert \rightarrow\infty}\left\vert u\left(
x,t\right)  -\left(  \sin\left(  \log\left\vert x\right\vert \right)
+\frac{1}{n}\cos\left(  \log\left\vert x\right\vert \right)  \right)
\right\vert =0. \label{u-conv-x}%
\end{equation}

\end{lemma}

\begin{remark}
Note that for fixed $t>0,\ $we have%
\[
\tilde{\alpha}:=\liminf_{\left\vert x\right\vert \rightarrow\infty}u\left(
x,t\right)  =\liminf_{\tau\rightarrow\infty}\varphi\left(  \tau\right)
=r=-\sqrt{1+\frac{1}{n^{2}}}<\liminf_{\tau\rightarrow\infty}H\left(
\tau\right)  =p=-1
\]
and%
\[
\tilde{\beta}:=\limsup_{\left\vert x\right\vert \rightarrow\infty}u\left(
x,t\right)  =\limsup_{\tau\rightarrow\infty}\varphi\left(  \tau\right)
=s=\sqrt{1+\frac{1}{n^{2}}}>\limsup_{\tau\rightarrow\infty}H\left(
\tau\right)  =q=1.
\]
Thus the\textbf{ space oscillation} of $u\left(  x,t\right)  \ $is, in
general, not necessarily bounded by the oscillation of $H\left(  \tau\right)
\ $(however, it is always bounded by the oscillation of $\varphi\left(
\tau\right)  $ due to the maximum principle).$\ $This is different\ from the
behavior for \textbf{time oscillation}, where we always have$\ r\leq
p\leq\alpha\leq\beta\leq q\leq s.$
\end{remark}

%

%TCIMACRO{\TeXButton{Proof}{\proof}}%
%BeginExpansion
\proof
%EndExpansion
By Remark \ref{rmk-rad}, we know that $u\left(  x,t\right)  \ $is\textbf{
radial} in $x\in\mathbb{R}^{n}\ $for each $t>0.\ $Therefore, it suffices to
look at the oscillation behavior of $u\left(  x,t\right)  $ for $\left\vert
x\right\vert \rightarrow\infty.\ $The representation formulas (\ref{rep}) and
(\ref{CONV-rad}) cannot be used here since they are valid only for $u\left(
0,t\right)  .$ Instead, we use the convolution formula (\ref{CONV}) to get
\begin{equation}
u\left(  x,t\right)  =\dfrac{1}{\left(  4\pi t\right)  ^{n/2}}\int
_{\mathbb{R}^{n}}e^{-\frac{\left\vert x-y\right\vert ^{2}}{4t}}\varphi\left(
\left\vert y\right\vert \right)  dy=\dfrac{1}{\pi^{n/2}}\int_{\mathbb{R}^{n}%
}e^{-\left\vert z\right\vert ^{2}}\varphi\left(  \left\vert x+\sqrt
{4t}z\right\vert \right)  dz,
\end{equation}
where by (\ref{phi-3}) we have$\ $
\[
\varphi\left(  \left\vert x+\sqrt{4t}z\right\vert \right)  =\left\{
\begin{array}
[c]{l}%
\sin\left(  \log\left(  \left\vert x+\sqrt{4t}z\right\vert +1\right)  \right)
%
%TCIMACRO{\TeXButton{vspace}{\vspace{3mm}}}%
%BeginExpansion
\vspace{3mm}%
%EndExpansion
\\
+\frac{\left\vert x+\sqrt{4t}z\right\vert }{n\left(  \left\vert x+\sqrt
{4t}z\right\vert +1\right)  }\cos\left(  \log\left(  \left\vert x+\sqrt
{4t}z\right\vert +1\right)  \right)  .
\end{array}
\right.
\]
By the Lebesgue Dominated Convergence Theorem and the identity%
\[
\sin\left(  \log\left(  \left\vert x+\sqrt{4t}z\right\vert +1\right)  \right)
=\left\{
\begin{array}
[c]{l}%
\sin\left(  \log\left\vert x\right\vert \right)  \cos\left(  \log\left(
\frac{\left\vert x+\sqrt{4t}z\right\vert +1}{\left\vert x\right\vert }\right)
\right)
%TCIMACRO{\TeXButton{vspace}{\vspace{3mm}}}%
%BeginExpansion
\vspace{3mm}%
%EndExpansion
\\
+\cos\left(  \log\left\vert x\right\vert \right)  \sin\left(  \log\left(
\frac{\left\vert x+\sqrt{4t}z\right\vert +1}{\left\vert x\right\vert }\right)
\right)  ,
\end{array}
\right.
\]
where
\[
\lim_{\left\vert x\right\vert \rightarrow\infty}\frac{\left\vert x+\sqrt
{4t}z\right\vert +1}{\left\vert x\right\vert }=1,\ \ \ \text{for
fixed\ }t\ \text{and }z\text{,}%
\]
we have%
\[
\lim_{\left\vert x\right\vert \rightarrow\infty}\left\vert \sin\left(
\log\left(  \left\vert x+\sqrt{4t}z\right\vert +1\right)  \right)
-\sin\left(  \log\left\vert x\right\vert \right)  \right\vert =0
\]
and so
\begin{equation}
\lim_{\left\vert x\right\vert \rightarrow\infty}\left\vert \dfrac{1}{\pi
^{n/2}}\int_{\mathbb{R}^{n}}e^{-\left\vert z\right\vert ^{2}}\sin\left(
\log\left(  \left\vert x+\sqrt{4t}z\right\vert +1\right)  \right)
dz-\sin\left(  \log\left\vert x\right\vert \right)  \right\vert =0.
\end{equation}
Similarly, we have%
\begin{equation}
\lim_{\left\vert x\right\vert \rightarrow\infty}\left\vert \dfrac{1}{\pi
^{n/2}}\int_{\mathbb{R}^{n}}e^{-\left\vert z\right\vert ^{2}}\left(
\tfrac{\left\vert x+\sqrt{4t}z\right\vert }{n\left(  \left\vert x+\sqrt
{4t}z\right\vert +1\right)  }\cos\left(  \log\left(  \left\vert x+\sqrt
{4t}z\right\vert +1\right)  \right)  \right)  dz-\frac{1}{n}\cos\left(
\log\left\vert x\right\vert \right)  \right\vert =0
\end{equation}
and (\ref{u-conv-x}) follows. The proof is done.$%
%TCIMACRO{\TeXButton{hfill}{\hfill}}%
%BeginExpansion
\hfill
%EndExpansion
\square$

\ \ \ \ \ \ 

\textbf{Acknowledgement.\ \ \ \ }Research supported by NCTS (National Center
for Theoretical Sciences)\ and MoST\ (Ministry of Science and Technology)\ of
Taiwan with grant number\ 108-2115-M-007-013-MY2.

\ \ \ \ \ \ \ \ 

\ \

\ \ \ \ \ \ \ \ \ 

\ \ \ \ \ \ \ \ \ \ \ 

Dong-Ho Tsai

Department of Mathematics

National Tsing Hua University

Hsinchu 30013,\ TAIWAN

E-mail:\ \textit{dhtsai@math.nthu.edu.tw} 

\ \ \ \ \ \ \ \ \ \ \


\begin{thebibliography}{99}                                                                                               %


\bibitem[CE]{CE}P. Collet,\ J.-P.\ Eckmann, \emph{Space-time behavior in
problems of hydrodynamic type:\ A case study},\ \textbf{Nonlinearity,
5\ }(1992) 1265-1302.

\bibitem[CT]{CT}M.-S.\ Chang,\ D.-H. Tsai, \emph{On the oscillation behavior
of solutions to the heat equation on }$\mathbb{R}^{n}$,\ \textbf{Journal
of\ Differential Equations, 268}\ (2020)\ 2040-2062.

\bibitem[E]{E}S.D. Eidel'man,\ \emph{Parabolic System}, North-Holland,
Amsterdam, 1969.

\bibitem[J]{J}F. John,\ \emph{Partial Differential Equations,\ 4th edition,}
\textbf{Applied Mathematical Sciences,\ v. 1},\ Springer-Verlag, 1982.

\bibitem[K]{K}S. Kamin, \emph{On stabilization of solutions of the Cauchy
problem for parabolic equations}, \textbf{Proc.\ Roy.\ Soc.\ Edinburgh Sect.
A, 76} (1976)\ 43-53.

\bibitem[N]{N}W.-M. Ni,\ \emph{The mathematics of diffusion,\ }%
\textbf{CBMS-NSF Regional Conference Series in Applied Math.,\ v. 82,\ SIAM} (2011).\ 

\bibitem[NT]{NT}M.\ Nara,\ M.\ Taniguchi,\ \emph{The condition on the
stability of stationary lines in a curvature flow in the whole plane}%
,\ \textbf{J. Diff. Eq. 237\ }(2007) 61-76.\ 

\bibitem[RE]{RE}V.D.\ Repnikov,\ S.D. Eidel'man, \emph{A new proof of the
theorem on the\ stabilization of the solution of the Cauchy problem for the
heat equation},\ \textbf{Math. USSR Sb.,\ 2} (1967)\ 135-139.

\bibitem[TN]{TN}D.-H. Tsai,\ C.-H.\ Nien, \emph{On the oscillation behavior of
solutions to the one-dimensional heat equation}, \textbf{Discrete \&
Continuous Dynamical Systems-A, vol.\ 39,\ no. 7} (2019)\ 4073-4089.
\end{thebibliography}
\end{document}